\documentclass[11pt]{article}
\usepackage{graphicx} 
\usepackage{amsfonts} 
\usepackage{amssymb}  
\usepackage{amsmath}
\usepackage{epstopdf}
\usepackage{hyperref}
\usepackage{amsthm}
\usepackage{color}

\begin{document}
\title{Isogeometric solution of Helmholtz equation with Dirichlet boundary condition: numerical experiences.}
\author{Victoria Hern\'andez Mederos $^a$, Isidro A. Abell\'o Ugalde $^b$,\\
 Rolando M. Bruno Alfonso $^c$, Domenico Lahaye $^d$, Valia Guerra Ones $^d$ \\ \\
$^a$ {\it Instituto de Cibern\'etica, Matem\'atica y F\'isica, ICIMAF, La Habana, Cuba}\\
$^b$ {\it CEPES, Universidad de La Habana, Cuba}\\
$^c$ {\it Facultad de Matem\'atica y Computaci\'on, Universidad de La Habana, Cuba}\\
$^d$ {\it DIAM, TU Delft, The Netherlands}
}
\date{}
\maketitle

\begin{abstract}
In this paper we use the Isogeometric method to solve the Helmholtz equation  with nonhomogeneous Dirichlet boundary condition over a bounded physical domain. Starting from the variational formulation of the problem, we show with details how to apply the isogeometric approach to obtain an approximation of the solution using biquadratic B-spline functions. To illustrate the power of the method we solve several difficult problems, which are particular cases of the Helmholtz equation, where the solution has discontinuous gradient in some points, or it is highly oscillatory. For these problems we explain how to select the knots of B-spline quadratic functions and how to insert knew knots in order to obtain good approximations of the exact solution on regions with irregular boundary. The results, obtained with a our Julia implementation of the method, prove that isogeometric approach produces approximations with a relative small error and computational cost.
\end{abstract}

{\bf keywords}: isogeometric analysis, numerical experiences, Helmholtz equation.

\section{Introduction}
In its most general form of Helmholtz equation is given by
\begin{equation}
-\triangle u(x,y) - k^2(x,y) u(x,y)=f(x,y),\;\;\;(x,y) \in \Omega \label{Heq}
\end{equation}
where $k(x,y)$ and $f(x,y)$ are known functions and $\triangle$ denotes the Laplacian.
The equation (\ref{Heq}) includes several interesting cases. For instance if $k\equiv 0$ then Helmholtz equation is reduced to Poisson equation. Moreover, the wave function that satisfies a Schr\"{o}dinger equation model of two interacting atoms \cite{Mit13} is solution of a Helmholtz equation with variable frequency $k(x,y)$.

Due to its importance in different fields such as acoustic, seismic and electromagnetic systems, the case $k$ equal to a positive constant have been intensively investigated over the years \cite{Ihl95}, \cite{Ihl97}, \cite{Ern11}, \cite{Moi14}. In this case, $u(x,y)$ is the amplitude of a wave traveling along $\Omega$ and $k$, called wavenumber, is the number of waves per unit of distance. When $k$ is small, the problem can be handled using low order Finite Element Method (FEM). But as the wavenumber is increased, low order finite elements become very expensive and many numerical difficulties appear.


Isogeometric analysis (IgA) was introduced by Hughes et al. in \cite{Hugh05} as an extension of FEM to solve partial differential equations (PDE). Its name reflects the fact that IgA uses B-splines functions for two different purposes: to parametrize the geometry $\Omega$ and to approximate the solution of the PDE. In comparison with FEM, IgA has two basic advantages: the boundary of the physical domain is represented exactly and the approximated solution of the PDE is smoother, with one or several continuous derivatives.

In this paper we apply the isogeometric approach to solve the Helmholtz equation (\ref{Heq}) over a bounded physical domain $\Omega$, with Lipschitz continuous boundary $\Gamma$ and nonhomogeneous Dirichlet boundary condition
\begin{equation}
u(x,y)=g(x,y), \;\;\; (x,y) \in \Gamma \label{Hebc}
\end{equation}
We are specially interested in problems where the physical regions $\Omega$ has irregular boundary, such as lagoons, lakes, islands, etc.
The first step to solve a PDE with IgA approach is the parametrization of $\Omega$ with a tensor product B-spline function.
This is currently an active research area, see for instance \cite{Xu10},
\cite{Xu11}, \cite{Ngu12}, \cite{Grav14}, \cite{Xu13a)}, \cite{Fall15},
\cite{Nian16}, \cite{Xu18} and \cite{Abe18}. In this paper, we parametrize  $\Omega$ with  a biquadratic tensor product B-spline, which is computed by the method proposed in \cite{Abe18}. This method produces high a quality parametrization of complex planar regions $\Omega$. This is very important because the quality of the parametrization has an strong influence on the speed of convergence of the approximated solution and also on the condition number of the discretization  matrix \cite{Xu11}, \cite{Nian16}.

The main contribution of this paper is the solution of three difficult problems obtained as particular cases of the Helmholtz equation using the IgA approach. In all cases the physical domain is a region with irregular boundary which increases the difficulty for computing the approximated solution. The first problem is a Poisson equation with a highly oscillatory solution. The second problem is also a Poisson equation, whose solution has discontinuous gradient in several points. The last problem is a Helmholtz equation  with variable frequency and a highly oscillatory solution near a singular point. In all the cases we show how to construct carefully the sequence of knots of the biquadratic B-spline functions, in order to obtain approximations to the exact solution having similar behavior, including singular points and oscillations. The performance of a Julia code to solve the previous problems in several regions with very irregular boundary are also included, showing that the method produces accurate approximations to the exact solution.

The paper is organized as follows. In section 2 we obtain the variational formulation of Helmholtz equation with nonhomogeneous boundary condition. The isogeometric method is described in section 3 using biquadratic B-spline functions to approximate the solution of the problem. Details about the basic steps of the method are given in this section, including the obtention of the linear system of equations, which provides the coefficients of the approximated solution written in the tensor product B-spline basis. Section 4 describes how to approximate the Dirichlet boundary condition. Moreover, computational aspects of the assembly process are also given.
In section 5 we show how to apply the IgA approach to the solution of several problems obtained from Helmholtz equation. Section 6 concludes the paper.

\section{Variational formulation. }
FEM and IgA have both the same theoretical basis, namely the weak or variational formulation of a PDE.
In this section we obtain the variational formulation of Helmholtz equation with homogeneous boundary condition.
Our problem with boundary condition (\ref{Hebc}) is reduced to a problem with homogeneous boundary condition
writing  the solution of (\ref{Heq}) as
\begin{equation}
u(x,y)=u_0(x,y)+u_g(x,y) \label{u}
\end{equation}
where the function $u_0$ satisfies (\ref{Heq}) and
\begin{equation}
u_0(x,y)=0, \;\; \mbox{ for} \;\; (x,y) \in \Gamma  \label{Hebcu0}
\end{equation}
while
\begin{equation}
u_g(x,y)=g(x,y),\;\; \mbox{ for} \;\;(x,y) \in \Gamma  \label{Hebcug}
\end{equation}
Thus, substituting (\ref{u}) in (\ref{Heq}) we transform the original problem (\ref{Heq})-(\ref{Hebc}) into the following
problem
\begin{equation}
-\triangle u_0(x,y) - k^2(x,y) u_0(x,y)=\widetilde{f}(x,y) \label{Hequ0}
\end{equation}
with homogeneous Dirichlet boundary condition, where

$$\widetilde{f}(x,y)=f(x,y)+\triangle u_g(x,y) + k^2(x,y)u_g(x,y).$$

Now, let $\mathcal{V}$ the Hilbert space of functions
\begin{equation}
\mathcal{V}=\{v \in H^1(\Omega),\;/v(x,y)=0\;\;\mbox{for}\;\;(x,y) \in \Gamma\} \label {V}
\end{equation}
which consists of all functions $v \in L_2(\Omega)$ that possess weak and square-integrable
first derivatives and that vanish on the boundary. The norm $\|v\|_{\mathcal{V}}$ in this space is given by
\begin{equation}
\|v\|_{\mathcal{V}}^2=\int\int_{\Omega} v^2 +\left(\frac{\partial v}{\partial x}\right)^2 +\left(\frac{\partial v}{\partial y}\right)^2 \; d\Omega \label{normV}
\end{equation}
To obtain the variational formulation we multiply (\ref{Hequ0}) by $v \in \mathcal{V}$ and integrate on $\Omega$
\begin{equation}
\int\int_{\Omega} (-\triangle u_0(x,y) - k^2(x,y) u_0(x,y))v(x,y)\;d\Omega= \int\int_{\Omega} \widetilde{f}(x,y)v(x,y)  \;d\Omega \label{eq0}
\end{equation}
Using in (\ref{eq0}) the Green formula
\begin{equation}
\int\int_{\Omega} \nabla u(x,y)^{t}\nabla v(x,y)\;d\Omega = -\int\int_{\Omega} \triangle u(x,y)v(x,y)\;d\Omega + \int_{\Gamma} \frac{\partial u}{\partial \mathbf{n}}v \; ds
\end{equation}
where $\nabla u =\left(\frac{\partial u}{\partial x},\frac{\partial u}{\partial y}\right)^{t}$ and $\mathbf{n}$ denotes the outer normal vector to $\Gamma$, we obtain
\begin{equation}
\int\int_{\Omega} (\nabla u_0(x,y)^{t}\nabla v(x,y)-k^2(x,y)u_0(x,y)v(x,y))\;d\Omega = \int\int_{\Omega} \widetilde{f}(x,y)v(x,y)  \;d\Omega + \int_{\Gamma} \frac{\partial u_0}{\partial \mathbf{n}}v(x,y) \; ds\label{eq1}
\end{equation}
Since $v \in \mathcal{V}$ the last integral in (\ref{eq1}) vanishes. Moreover, we can use the Green formula again to simplify the right hand side in (\ref{eq1}) obtaining
\begin{equation}
\int\int_{\Omega} \widetilde{f}(x,y)v(x,y)  \;d\Omega= \int\int_{\Omega} (f(x,y) +k^2(x,y)u_g(x,y))v(x,y) \;d\Omega
- \int\int_{\Omega} \nabla u_g(x,y)^{t}\nabla v(x,y)  \;d\Omega
\label{eq2}
\end{equation}
Finally, substituting (\ref{eq2}) in (\ref{eq1}) we obtain the variational formulation: find $u_0 \in \mathcal{V}$ such that for all $v \in \mathcal{V}$
\begin{equation}
a(u_0,v)=G(v)
\label{variational}
\end{equation}
where $a(u,v)$ is the bilinear form
\begin{equation}
a(u,v)=\int\int_{\Omega} (\nabla u(x,y)^{t}\nabla v(x,y)-k^2(x,y)u(x,y)v(x,y))\;d\Omega
\label{auv}
\end{equation}
and $G(v)$ is the linear form
\begin{equation}
G(v)=\int\int_{\Omega} (f(x,y) +k^2(x,y)u_g(x,y))v(x,y) \;d\Omega
- \int\int_{\Omega} \nabla u_g(x,y)^{t}\nabla v(x,y)  \;d\Omega
\label{Fv}
\end{equation}


The existence and uniqueness of weak solution has been very well studied when $k(x,y)^2=\lambda$, where $\lambda$ is a real constant, see for instance \cite{Spen15}. For $\lambda=0$ the bilinear form $a(u,v)$ given by (\ref{auv}) is {\it coercive}. Therefore, Lax-Milgram theorem guarantees the existence and uniqueness of a solution to the variational problem (\ref{variational}) and continuous dependence of the solution on the data. On the other hand, if $\lambda = \lambda_j$, where $\lambda_j$ is the $j$-th Dirichlet eigenvalue of the negative Laplacian in $\Omega$, i.e. there exists a $u_j \in H^1(\Omega)\setminus 0$ such that $-\triangle u_j = \lambda_j u_j$ in $\Omega$ and $u_j=g$ on $\Gamma $, then the problem has solution but it is not unique. Finally, if $\lambda$ is not an eigenvalue of the negative Laplacian, then the bilinear form $a(u,v)$ satisfies a G\"{a}rding inequality and again the variational problem (\ref{variational}) has a unique solution which depends continuously on $f$.

There are few results in the literature about the Helmholtz equation with variable coefficient $k(x,y)$. In the recent paper \cite{Gra18}, existence and uniqueness results for this problem are obtained under rather general conditions on the function $k(x,y)$, using the unique continuation principle and the Fredholm alternative.

\section{Galerkin method with isogeometric approach.}

The Galerkin method replaces the infinite-dimensional space $\mathcal{V}$ by a finite-dimensional subspace $\mathcal{V}_h \subset \mathcal{V}$ and solves the corresponding discrete problem. In the classical FEM, the subspace $\mathcal{V}_h$ consists of piecewise polynomials with global $C^0$ continuity. This space is defined in terms of a partition of the physical domain
$\Omega$ in a mesh of triangles or quadrilaterals. In the isogeometric approach \cite{Cott09},  the subspace $\mathcal{V}_h$ is generated by tensor product B-spline functions ( or more general by NURBs functions ) with higher global continuity. Moreover, it is assumed that the physical domain $\Omega$ is topologically equivalent to the unit square $\hat{\Omega}$, thus its boundary can be divided into 4 curves in such a way that consecutive curves are the image by a parametrization
$$\mathbf{F}(\xi,\eta):\hat{\Omega}\longrightarrow\Omega$$
of consecutive sides of $\hat{\Omega}$. In this paper, we assume that $\mathbf{F}(\xi,\eta)$ is an injective biquadratic B-spline function that can be written as \cite{Boor01}
\begin{equation}
\label{MapBiCuadExpli}
\mathbf{F}(\xi,\eta)=(x(\xi,\eta),y(\xi,\eta))^t=\sum_{i=1}^n\sum_{j=1}^m \mathbf{P}_{i,j}B_{i,t^\xi}^3(\xi)B_{j,t^\eta}^3(\eta)
\end{equation}
where $\mathbf{P}_{i,j}=(P_{i,j}^x,P_{i,j}^y)^{t},\;i=1,...,n,\,j=1,...,m$ are the control points, $B_{i,t^\xi}^3(\xi)$ is the $i$-th quadratic B-spline for the knot sequence $t^{\xi}$ and  $B_{j,t^\eta}^3(\eta)$  is the $j$-th  quadratic B-spline for the knot sequence $t^{\eta}$ with
\begin{eqnarray}
t^{\xi}&=&(0,0,\xi_1,\xi_2,...,\xi_{n-1},1,1),\;\;\;0=\xi_1<\xi_2<...<\xi_{n-1}=1 \label{tchi_cuad} \\
t^{\eta}&=&(0,0,\eta_1,\eta_2,...,\eta_{m-1},1,1),\;\;\;0=\eta_1<\eta_2<...<\eta_{m-1}=1 \label{teta_cuad}
\end{eqnarray}
In other words, $\mathbf{F}(\xi,\eta)$ is a function in $\mathbb{S}_{3,t^\xi} \bigotimes \mathbb{S}_{3,t^\eta}$, where $\mathbb{S}_{3,t}$ denotes the space of quadratic spline functions for the knot sequence $t$.
To simplify the notation, in the rest of the paper we don't write the subindex $t^\xi$ or $t^\eta$ of the B-spline functions. The functions
\begin{equation}
B_{i,j}^3(\xi,\eta):=B_i^3(\xi)B_j^3(\eta),\;\;i=1,...n,\;\;j=1,...,m
\end{equation}
define a basis of $\mathbb{S}_{3,t^\xi} \bigotimes \mathbb{S}_{3,t^\eta}$.
Then, due to the assumptions on the parameterization $\mathbf{F}$, the functions
\begin{equation}
\Phi_{i,j}(x,y):=(B_{i,j}^3\, o\, \mathbf{F}^{-1})(x,y),\;\;i=1,...n,\;\;j=1,...,m
\label{Phi}
\end{equation}
are independent in $\Omega$. The control points  $\mathbf{P}_{i,j},\;i=1,...,n,\,j=1,...,m$ are computed as the vertices of a quadrilateral mesh which is obtained by minimizing a functional \cite{Abe18}.
\smallskip
With the help of $\mathbf{F}$, integrals (\ref{auv}),(\ref{Fv}) over $\Omega$ can be transformed into integrals over $\hat{\Omega}$ by means of the integration rule
$$\int \int_{\Omega} h(x,y)\;d\Omega=\int_0^1 \int_0^1 h(\mathbf{F}(\xi,\eta)) |\det J\mathbf{F}(\xi,\eta) |\;\;d \xi\,d \eta $$
where $J\mathbf{F}$ denotes the Jacobian matrix of the parametrization,
\begin{equation}
J\mathbf{F}(\xi,\eta)=
\left(
     \begin{array}{cc}
       x_{\xi} & x_{\eta} \\
       y_{\xi} & y_{\eta} \\
     \end{array}
   \right)
\label{Jacobian}
\end{equation}

Applying the chain rule to $h(x,y)= h(\mathbf{F}(\xi,\eta))$ we obtain
$$\nabla_{(x,y)}h(x, y) = J\mathbf{F}(\xi, \eta)^{-t}\nabla_{(\xi,\eta)}h(\xi,\eta)$$
where the notation $\nabla_{(x,y)}$ means that partial derivatives are computed with respect to variables $x,y$.
Hence, integrals of the weak form (\ref{auv}),(\ref{Fv}) can be written as
\begin{eqnarray}
a(u,v)&=&\int_0^1\int_0^1(J\mathbf{F}(\xi, \eta)^{-t}\nabla_{(\xi,\eta)}u)^{t}
                         (J\mathbf{F}(\xi, \eta)^{-t}\nabla_{(\xi,\eta)}v) \;|\det J\mathbf{F}(\xi,\eta) |\;d \xi d \eta \nonumber\\
     &-&\int_0^1\int_0^1 k^2(\mathbf{F}(\xi, \eta))u(\mathbf{F}(\xi, \eta))v(\mathbf{F}(\xi, \eta))
                         |\det J\mathbf{F}(\xi,\eta) |\;d \xi\,d \eta
\label{axieta}
\end{eqnarray}
and
\begin{eqnarray}
G(v)&=&\int_0^1\int_0^1 (f(\mathbf{F}(\xi,\eta)) +k^2(\mathbf{F}(\xi,\eta))u_g(\mathbf{F}(\xi,\eta)))v(\mathbf{F}(\xi,\eta))
\;|\det J\mathbf{F}(\xi,\eta) |\;d \xi d \eta \nonumber\\
&-& \int_0^1\int_0^1 (J\mathbf{F}(\xi, \eta)^{-t}\nabla_{(\xi,\eta)}u_g(\xi,\eta))^{t}
J\mathbf{F}(\xi, \eta)^{-t}\nabla_{(\xi,\eta)}v(\xi,\eta)
\;|\det J\mathbf{F}(\xi,\eta) |\;d \xi d \eta
\label{Fxieta}
\end{eqnarray}

The approximate solution $u_0^h(x,y)$ is sought in the space
\begin{equation}
\mathcal{V}_h=\{span(\Phi_{i,j}(x,y))\,\;\mbox{ such that}\;\; \Phi_{i,j}(x,y)=0,\;\;\mbox{for} \;\; (x,y) \in \Gamma\}
\label{Vh}
\end{equation}

Taking into account (\ref{tchi_cuad}) and (\ref{teta_cuad}) it is easy to check that \cite{Boor01}
\begin{equation}
\Phi_{i,j}(x,y)=0,\;\;\;\;(x,y)\in \Gamma, \;\;\mbox{for}\;\;i=2,...,n-1,\;j=2,...,m-1
\label{baseVh}
\end{equation}
From (\ref{Vh}) and (\ref{baseVh}) we conclude that
\begin{equation}
\mathcal{V}_h=span\{\Phi_{i,j}(x,y)\},\;\;\;\mbox{for}\;\;i=2,...,n-1,\;j=2,...,m-1
\label{VhPhi}
\end{equation}
Hence, $u_0^h(x,y)$ can be written as
\begin{equation}
u_0^h(x,y)=\sum_{i=1}^{n}\sum_{j=1}^{m}\gamma_{i,j}\Phi_{i,j}(x,y)
\label{u0comp}
\end{equation}
where $\gamma_{i,j} \in \mathbb{R},\;i=1,..,n,\;j=1,...,m$ and the coefficients $\gamma_{1,j},\gamma_{n,j},\;j=1,...,m$ and $\gamma_{i,1},\gamma_{i,m},\;i=1,...,n$ must be forced to be zero. In order to obtain a linear system for the unknowns $\gamma_{i,j}$ it is convenient to vectorize the basis functions and the corresponding coefficients in (\ref{u0comp}) introducing the change of index

\begin{equation}
p:=n(j-1)+i,\;\;i=1,...,n,\;\;j=1,...,m
\label{indexnew}
\end{equation}
With this transformation, the expression (\ref{u0comp}) can be written as
\begin{equation}
u_0^h(x,y)=\sum_{p=1}^{N}\alpha_p \psi_p(x,y)
\label{u0pnew}
\end{equation}
where $N:=nm$ and
\begin{eqnarray}
\mathbf{\alpha}:&=&(\alpha_1,...,\alpha_N)^{t}=(\gamma_{1,1},...,\gamma_{n,1},...,
\gamma_{1,m},...,\gamma_{n,m})^{t}  \label{alpha}\\
\mathbf{\psi}(x,y):&=&(\psi_1(x,y),...,\psi_N(x,y))\\
 &=&(\phi_{1,1}(x,y),...,\phi_{n,1}(x,y),...,\phi_{1,m}(x,y),...,\phi_{n,m}(x,y)) \label{psi}
\end{eqnarray}
\smallskip
We subdivide the set of indexes $I=\{1,2,...,N\}$ in two subsets:
$I=I_0 \cup I_1$, where $I_0$ is the set of basis functions that generate $\mathcal{V}_h$ and $I_1$ is the set of basis functions that are different from $0$ on $\Gamma$. In other words, $I_0$ is the set of the indexes
(\ref{indexnew}) corresponding to basic functions of $\mathcal{V}_h$: $\phi_{i,j}(x,y),\;i=2,...,n-1,$ $\;j=2,...,m-1$ or $\mathcal{V}_h=span\{\psi_p(x,y),\;p \in I_0\}$. Moreover, $I_1$ is the set of indexes (\ref{indexnew}) corresponding to functions $\phi_{1,j}(x,y),\phi_{n,j}(x,y),\;j=1,...,m$ and $\phi_{i,1}(x,y),\phi_{i,m}(x,y),\;i=1,...,n$.
\smallskip

Substituting in (\ref{variational}), the expressions (\ref{axieta}) and (\ref{Fxieta})
and also $u(x,y)$ by $u_0^h(x,y)$ given by (\ref{u0pnew}) and $v(x,y)$ by the basis function of $\mathcal{V}_h$, $\psi_q(x,y),\;q \in I_0$, we obtain the Galerkin formulation: find $u_0^h(x,y)$ such that

\begin{eqnarray*}
&& \int_0^1\int_0^1 \left[\left(\sum_{p=1}^N \alpha_p \nabla \psi_p\right)^{t}(J\mathbf{F}^{t}J\mathbf{F})^{-1} \nabla \psi_q -k^2 \left(\sum_{p=1}^N \alpha_p \psi_p\right)\psi_q\right]\,|\det J\mathbf{F}|\;d \xi d \eta\\
=& &\int_0^1\int_0^1 \left[(f+k^2u_g)\psi_q - (\nabla u_g)^t(J\mathbf{F}^{t}J\mathbf{F})^{-1} \nabla \psi_q\right]\,|\det J\mathbf{F}|\;d \xi d \eta, \qquad q \in I_0
\end{eqnarray*}
where we have simplified the expressions omitting the dependence of $(\xi,\eta)$ of all
functions. The last expression is equivalent to
\begin{eqnarray*}
& & \sum_{p=1}^N\left[ \int_0^1\int_0^1 \left[ \left(\nabla \psi_p\right)^{t}(J\mathbf{F}^{t}J\mathbf{F})^{-1} \nabla \psi_q -k^2 \psi_p\psi_q\right]\,|\det J\mathbf{F}|\;d \xi d \eta \right] \alpha_p \\
=& &\int_0^1\int_0^1 \left[(f+k^2u_g)\psi_q - (\nabla u_g)^t(J\mathbf{F}^{t}J\mathbf{F})^{-1} \nabla \psi_q\right]\,|\det J\mathbf{F}|\;d \xi d \eta,\;\;q \in I_0
\end{eqnarray*}
These equations can be written in matrix form as
\begin{equation}
\mathbf{A}\mathbf{\alpha}=\mathbf{b}
\label{sist}
\end{equation}
where
\begin{equation}
\mathbf{A}=(a_{q,p})=\int_0^1\int_0^1 \left[ \left(\nabla \psi_p\right)^{t}(J\mathbf{F}^{t}J\mathbf{F})^{-1} \nabla \psi_q -k^2 \psi_p\psi_q\right]\,|\det J\mathbf{F}|\;d \xi d \eta, \;\; p=1,...,N,\;q \in I_0
\label{A}
\end{equation}
\begin{equation}
\mathbf{b}=(b_q)=\int_0^1\int_0^1 \left[(f+k^2u_g)\psi_q - (\nabla u_g)^t(J\mathbf{F}^{t}J\mathbf{F})^{-1} \nabla \psi_p\right]\,|\det J\mathbf{F}|\;d \xi d \eta,\;q \in I_0
\label{b}
\end{equation}
and $\mathbf{\alpha}$ is given by (\ref{alpha}). The unknown coefficients $\alpha_1,...,\alpha_N$ are computed solving the linear system $\widetilde{\mathbf{A}}\mathbf{\alpha}=\widetilde{\mathbf{b}}$, where the rows of $\widetilde{\mathbf{A}}=\widetilde{a}_{q,p}$ and $\widetilde{\mathbf{b}}=(\widetilde{b}_q)$ corresponding to basic functions of $\mathcal{V}_h$ ( $\psi_q,\;q \in I_0$) are given by (\ref{A})-(\ref{b}). To guarantee that the coefficients $\alpha_q,\;q \in I_1$ are zero (i.e the coefficient of functions $\psi_q,\;q \in I_1$) we set
$\widetilde{a}_{q,q}=1, \widetilde{a}_{q,p}=0,\;p \neq q, \;q\in I_1$  and $\widetilde{b}_q=0,\;q \in I_1$.

\section{Computing the B-spline approximated solution.}

In this section we explain how to compute a B-spline approximation of the function $u_g(x,y)$. Moreover, we give some details about the efficient implementation of the procedure to compute the global matrix and the right hand side vector of the linear system (\ref{sist}) whose solution are the B-spline coefficients of the approximated solution $u_0^h(x,y)$.

\subsection{Approximating the function $u_g(x,y)$.}\label{seccionug}

The function $u_g(x,y)$ satisfying boundary condition (\ref{Hebcug}) is approximated by a function $u_g^h(x,y)$ in  $\mathbb{S}_{3,t^\xi} \bigotimes \mathbb{S}_{3,t^\eta}$ written as
\begin{equation}
u_g^h(x,y)=\sum_{i=1}^{n}\sum_{j=1}^{m}\delta_{i,j}\Phi_{i,j}(x,y)
\label{ugcomp}
\end{equation}
with $\Phi_{i,j}(x,y)$  given by (\ref{Phi}). The unknown coefficients $\delta_{i,j},\;i=1,\cdots, n,\;j=1, \cdots, m$ are computed requiring that $u_g^h(x,y)$ interpolates the function $g(x,y)$, defining the Dirichlet boundary condition, at a sequence of points on $\Gamma$. More precisely, we select as interpolation sites $\widetilde{\xi}_i$ and $\widetilde{\eta}_j$ in the directions $\xi$ and $\eta$ respectively, the Greville abscissas, which in this case are the averages of 2 successive knots in the sequences $t^{\xi}$ and $t^{\eta}$:
\begin{eqnarray}
\widetilde{\xi}_k&=&\frac{t^{\xi}_{k+1}+t^{\xi}_{k+2}}{2},\;\;\;\;k=1, \cdots, n \label{intpxi}\\
\widetilde{\eta}_l&=&\frac{t^{\eta}_{l+1}+t^{\eta}_{l+2}}{2},\;\;\;\;l=1, \cdots, m \label{intpeta}
\end{eqnarray}
Evaluating the map $\mathbf{F}(\xi,\eta)$ given by (\ref{MapBiCuadExpli}) we obtain the sequence of interpolating points on $\Gamma$:
\begin{eqnarray*}
(x_k^{\it s},y_k^{\it s})&:=&\mathbf{F}(\widetilde{\xi}_k,0),\;\;\;\;(x_k^{\it n},y_k^{\it n}):=\mathbf{F}(\widetilde{\xi}_k,1),\;\;\;k=1,\cdots,n\\
(x_l^{\it w},y_l^{\it w})&:=&\mathbf{F}(0,\widetilde{\eta}_l),\;\;\;(x_l^{\it e},y_l^{\it e}):=\mathbf{F}(1,\widetilde{\eta}_l),\;\;l=1,\cdots,m
\end{eqnarray*}


Observe that points $(x_k^{\it s},y_k^{\it s}),(x_k^{\it n},y_k^{\it n}),$ $\;k=1,\cdots,n$ are on the ``south" and ``north" boundaries of $\Omega$, i.e in the boundary curves $\mathbf{F}(\xi,0)$ and $\mathbf{F}(\xi,1)$ respectively. Similarly, points $(x_l^{\it w},y_l^{\it w}),(x_l^{\it e},y_l^{\it e})$ $\;l=1,\cdots,m$ are on the ``west" and ``east" boundaries of $\Omega$, i.e in the boundary curves $\mathbf{F}(0,\eta)$ and $\mathbf{F}(1,\eta)$ respectively. In consequence, boundary coefficients of $u_g^h(x,y)$ in (\ref{ugcomp}) are computed from the interpolation conditions:
\begin{eqnarray*}
u_g^h(x_k^s,y_k^s)&=&g(x_k^s,y_k^s),\;\;\;\;\;u_g^h(x_k^n,y_k^n)=g(x_k^n,y_k^n),\;\;\;k=1,\cdots,n\\
u_g^h(x_l^w,y_l^w)&=&g(x_l^w,y_l^w),\;\;\;\;\;u_g^h(x_l^e,y_l^e)=g(x_l^e,y_l^e),\;\;\;l=1,\cdots,m
\end{eqnarray*}
Taking into account that the boundary knots in the sequences (\ref{tchi_cuad}) and (\ref{teta_cuad}) have multiplicity 3, from (\ref{ugcomp}) and (\ref{Phi}) we obtain that the previous interpolation conditions
can be written as,
\begin{eqnarray}
\sum_{i=1}^{n}\delta_{i,1} B_i^3(\widetilde{\xi}_k)&=&g(x_k^s,y_k^s),\;\;\;\;\;\;
\sum_{i=1}^{n}\delta_{i,m} B_i^3(\widetilde{\xi}_k)=g(x_k^n,y_k^n), \;\;k=1,\cdots,n \label{sys1}\\
\sum_{j=1}^{m}\delta_{1,j}B_j^3(\widetilde{\eta}_l)&=&g(x_l^w,y_l^w),\;\;\;\;\;\;
\sum_{j=1}^{m}\delta_{n,j}B_j^3(\widetilde{\eta}_l)=g(x_l^e,y_l^e),\;\;l=1,\cdots,m \label{sys3}
\end{eqnarray}

Observe that the matrix $B_1:=[B_i^3(\widetilde{\xi}_k)]_{i,k=1}^n$ of linear systems (\ref{sys1}) is the same and also the linear systems (\ref{sys3}) have the same matrix $B_2:=[B_j^3(\widetilde{\eta}_l)]_{j,l=1}^m$. Matrices $B_1$ and $B_2$ are nonsingular since hypothesis of Shoenberg-Whitney theorem \cite{Boor01} hold for interpolation sites (\ref{intpxi}) and (\ref{intpeta}). Thus, coefficients $\delta_{i,1},\delta_{i,m},\;i=1,\cdots,n$ and $\delta_{1,j},\delta_{n,j},\;j=1,\cdots,m $ are computed solving the corresponding linear systems.
The rest of the coefficients $\delta_{i,j},\;i=2,\cdots,n-1,\;j=2,\cdots,m-1$ are assigned as zero.

\begin{figure}[htb]
\center
\includegraphics[scale=0.3]{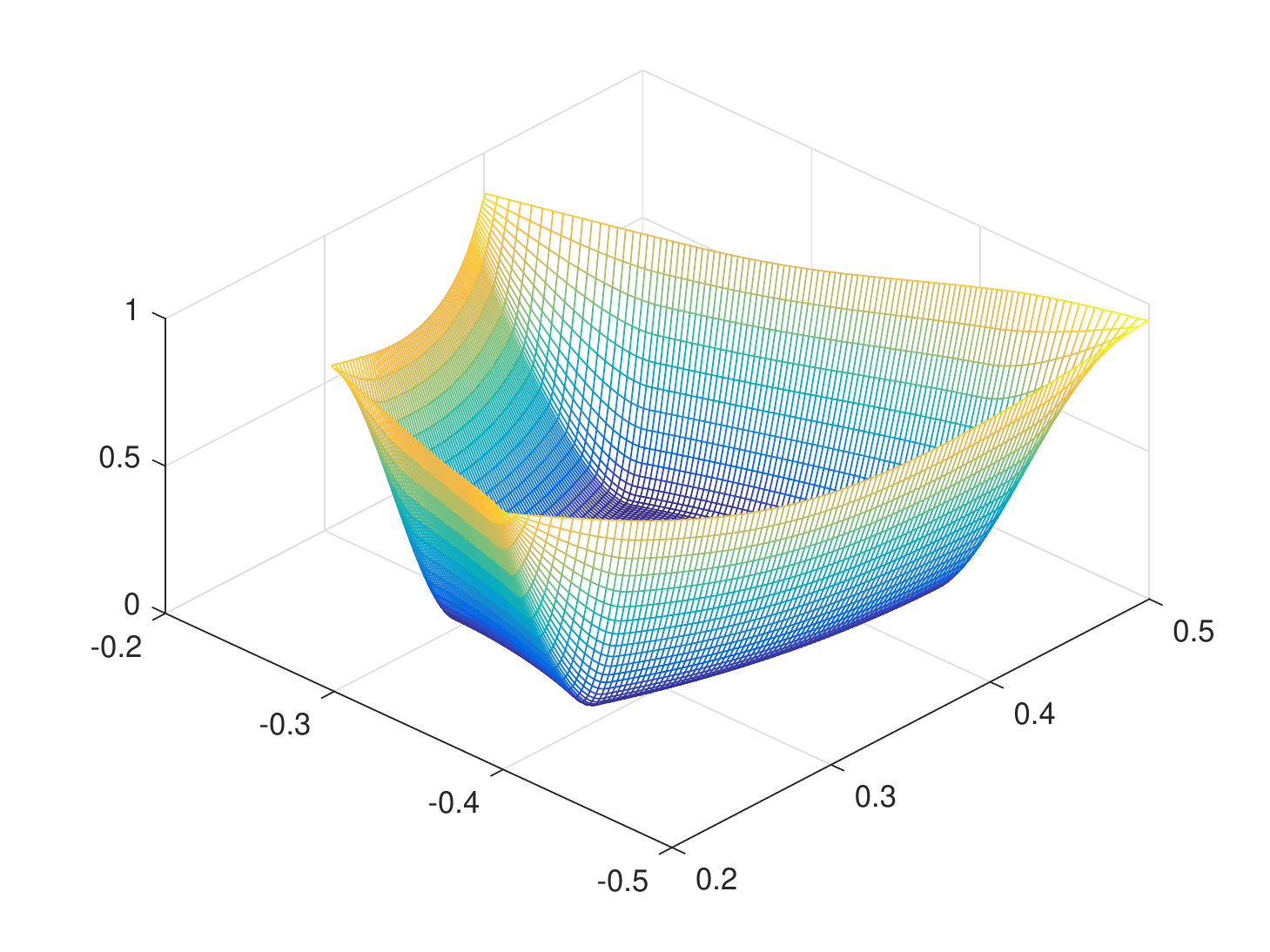}

\caption{Function $u_g^h(x,y)$ satisfying Dirichlet boundary condition for the region considered in the first example of table 1.}
\label{Fig:ug}
\end{figure}

Figure \ref{Fig:ug} shows the graphics of function $u_g^h(x,y)$ for the region considered in the first example of table 1.

\subsection{Assembling the global stiffness matrix and the right-hand side vector.}

The process of building the global stiffness matrix (\ref{A}) and the force vector (\ref{b}) is known in the FEM literature as {\it assembly}. This process doesn't not compute the elements of matrix $\mathbf{A}$ and vector $\mathbf{b}$, one entry at a time, as a first glance at the formulation (\ref{A})-(\ref{b}) might imply. Instead, one loops through the elements $\widehat{\Omega}_{k,l}:=[\xi_k,\xi_{k+1}] \times [\eta_l,\eta_{l+1}]$, building as we go local stiffness matrices and vectors $\mathbf{A}^{k,l}$ and $\mathbf{b}^{k,l}$ respectively, for $k=1,...,n-2,\;l=1,...,m-2$. Every entry of each of these dense matrices and vectors is then added to the appropriate spot in the global stiffness matrix $\mathbf{A}$ and vector $\mathbf{b}$.

Since in our problem the basic functions are biquadratic B-splines, only $9$ basic functions are different from zero in $\widehat{\Omega}_{k,l}$. These functions are

\begin{equation}
(\Phi_{k,l},\Phi_{k,l+1},\Phi_{k,l+2},\Phi_{k+1,l},\Phi_{k+1,l+1},\Phi_{k+1,l+2},
\Phi_{k+2,l},\Phi_{k+2,l+1},\Phi_{k+2,l+2})
\label{baselocal}
\end{equation}

Therefore, each local matrix $\mathbf{A}^{k,l}$ and the corresponding vector $\mathbf{b}^{k,l}$ are of order $9 \times 9$ and $9$ respectively. Denote by $p_1,...,p_9$ the global index of basic functions (\ref{baselocal}) computed using (\ref{indexnew}). Then
\begin{equation}
\mathbf{A}^{k,l}=
\left(
  \begin{array}{ccc}
    I_A(\psi_{p_1},\psi_{p_1}) & \cdots & I_A(\psi_{p_1},\psi_{p_9}) \\
    I_A(\psi_{p_2},\psi_{p_1}) & \cdots & I_A(\psi_{p_2},\psi_{p_9}) \\
    \vdots & \vdots & \vdots \\
    I_A(\psi_{p_9},\psi_{p_1}) & \cdots & I_A(\psi_{p_9},\psi_{p_9}) \\
  \end{array}
\right)
\end{equation}
where for $i,j=1,...,9$
\begin{equation}
I_A(\psi_{p_i},\psi_{p_j})=
\int_{\xi_k}^{\xi_{k+1}}\int_{\eta_l}^{\eta_{l+1}} \left[ \left(\nabla \psi_{p_i}\right)^{t}(J\mathbf{F}^{t}J\mathbf{F})^{-1} \nabla \psi_{p_j} -k^2 \psi_{p_i}\psi_{p_j}\right]\,|\det J\mathbf{F}|\;d \xi d \eta
\label{IA}
\end{equation}
Similarly,
\begin{equation}
\mathbf{b}^{k,l}=\left(
                   \begin{array}{c}
                      I_b(\psi_{p_1}) \\
                     I_b(\psi_{p_2}) \\
                     \vdots \\
                     I_b(\psi_{p_9}) \\
                   \end{array}
                 \right)
\end{equation}
where for $i=1,...,9$
\begin{equation}
I_b(\psi_{p_i})=\int_{\xi_k}^{\xi_{k+1}}\int_{\eta_l}^{\eta_{l+1}} \left[(f+k^2u_g)\psi_{p_i} - (\nabla u_g)^t(J\mathbf{F}^{t}J\mathbf{F})^{-1} \nabla \psi_{p_i}\right]\,|\det J\mathbf{F}|\;d \xi d \eta
\label{Ib}
\end{equation}

The integrals (\ref{IA}) and (\ref{Ib}) are computed approximately using Gaussian quadratures \cite{Hug10}.
Observe that $I_A(\psi_{p_i},\psi_{p_j}),\;i,j=1,...,9$ must be added up in $\mathbf{\widetilde{A}}(p_i,p_j)$. Similarly, $I_b(\psi_{p_i}),\;i=1,...,9$ must be added up $\mathbf{\widetilde{b}}(p_i)$.

\smallskip
Finally, the approximated solution $u^h(x,y)$ of the problem is given by $u^h(x,y)=u_0^h(x,y)+u_g^h(x,y)$. From (\ref{ugcomp}) and (\ref{u0comp}) it follows that
\begin{equation}
u^h(x,y)=\sum_{i=1}^{n}\sum_{j=1}^{m}\beta_{i,j}\Phi_{i,j}(x,y)
\label{uhcomp}
\end{equation}
where $\beta_{i,j}=\delta_{i,j}+\gamma_{i,j},\;i=1,\cdots, n,\;j=1, \cdots, m$.

\section{Numerical results}

In this section we describe our experiences solving the Helmholtz equation with Dirichlet boundary conditions using IgA approach. Our study includes the Poisson equation and the Helmholtz equation with variable frequency. In all the cases, the exact solution is known and therefore it is possible to compute the numerical error.  In our experiments we compute the $L_2\;error$ of the approximated solution $u^h(x,y)$ given by
\begin{equation}
(L_2\;error)^2=\int_0^1\int_0^1 (u(\mathbf{F}(\xi,\eta)) - u^h(\mathbf{F}(\xi,\eta)))^2 d\xi\,d\eta
\label{l2error}
\end{equation}
and also the $H_1\;error$ in the norm (\ref{normV}) given by
\begin{eqnarray}
(H_1\;error)^2&=&\int_0^1\int_0^1 (u(\mathbf{F}(\xi,\eta)) - u^h(\mathbf{F}(\xi,\eta)))^2 d\xi\,d\eta  + \int_0^1\int_0^1 \left(\frac{\partial u(\mathbf{F}(\xi,\eta))}{\partial \xi}- \frac{\partial u^h(\mathbf{F}(\xi,\eta))}{\partial \xi}\right)^2 d\xi\,d\eta \nonumber \\
&+& \int_0^1\int_0^1 \left(\frac{\partial u(\mathbf{F}(\xi,\eta))}{\partial \eta}-\frac{\partial u^h(\mathbf{F}(\xi,\eta))}{\partial \eta} \right)^2 d\xi\,d\eta
\label{H1error}
\end{eqnarray}

We consider several physical domains, with emphasis in planar regions $\Omega$ with irregular boundaries. The numerical results reported here have been obtained with the help of our computational implementation of isogeometric method in {\it Julia} language. This implementation uses biquadratic B-splines functions and computes the control points of  the map $\mathbf{F}(\xi,\eta)$ that parametrizes $\Omega$  by minimizing a functional \cite{Abe18}. We run our experiments in a PC with i5 processor and 4Gb of RAM.

\subsection{Oscillatory Poisson equation}

Our first example is the Poisson equation:
$$ -\triangle u(x,y) = 2\pi^2 \sin(\pi x) \sin(\pi y), \;\;(x,y) \in \Omega$$
The exact solution of this problem is
\begin{equation}
u(x,y) = \sin(\pi x) \sin(\pi y)
\label{uexacsin}
\end{equation}
a function that is highly oscillatory in $\Omega$. The function $g(x,y)$ defining the Dirichlet boundary condition (\ref{Hebc}) is $g(x,y)=u(x,y),\;\;(x,y) \in \Gamma$. The vector field of function (\ref{uexacsin}) is given by
$$ \nabla u =\left(\frac{\partial u}{\partial x},\frac{\partial u}{\partial y}\right)^{t}=(\pi \cos(\pi x) \sin (\pi y),\pi \sin(\pi x) \cos(\pi y))^{t}$$

\begin{figure}[htb]
\center
\includegraphics[scale=0.3]{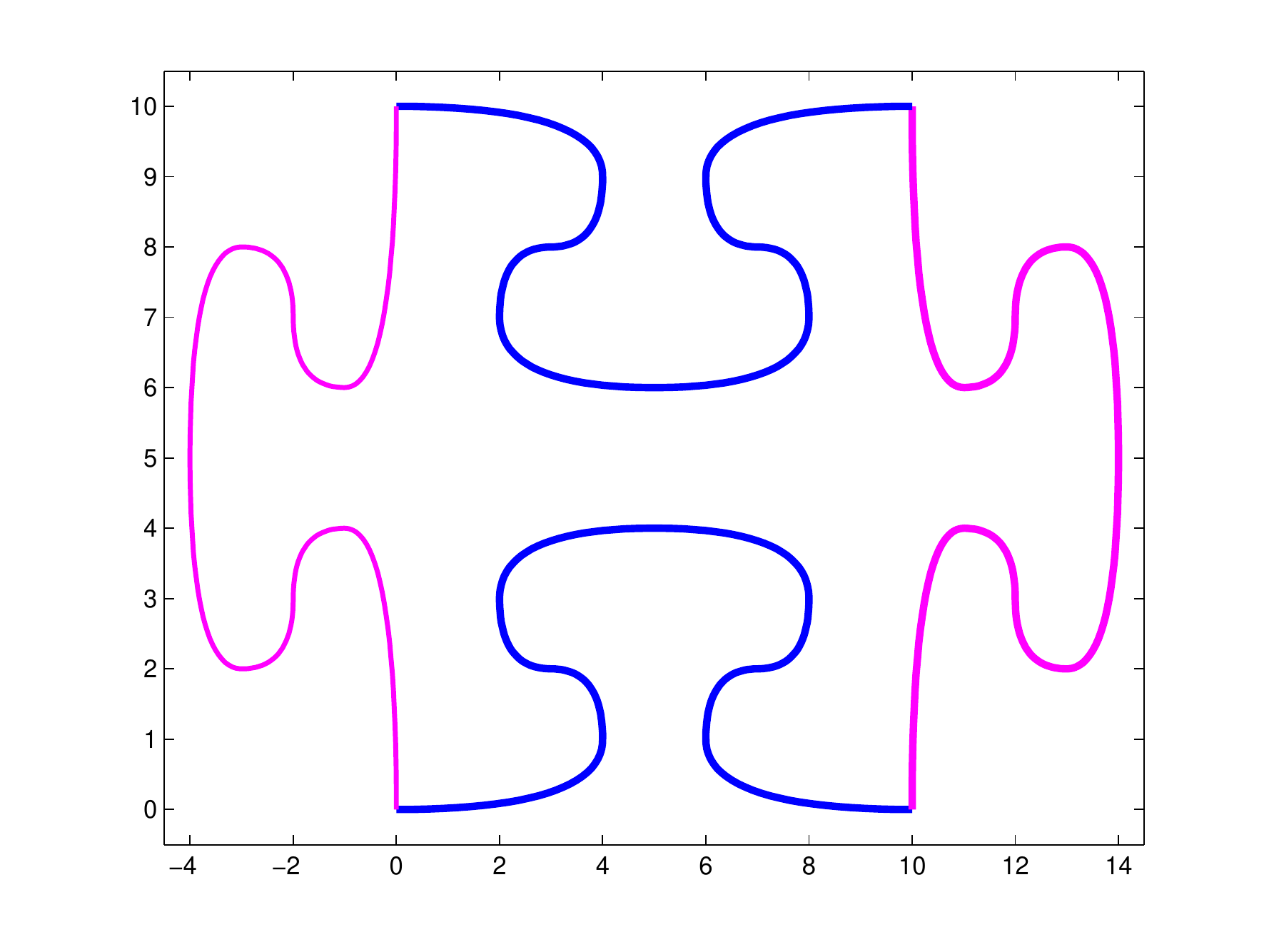}
\includegraphics[scale=0.3]{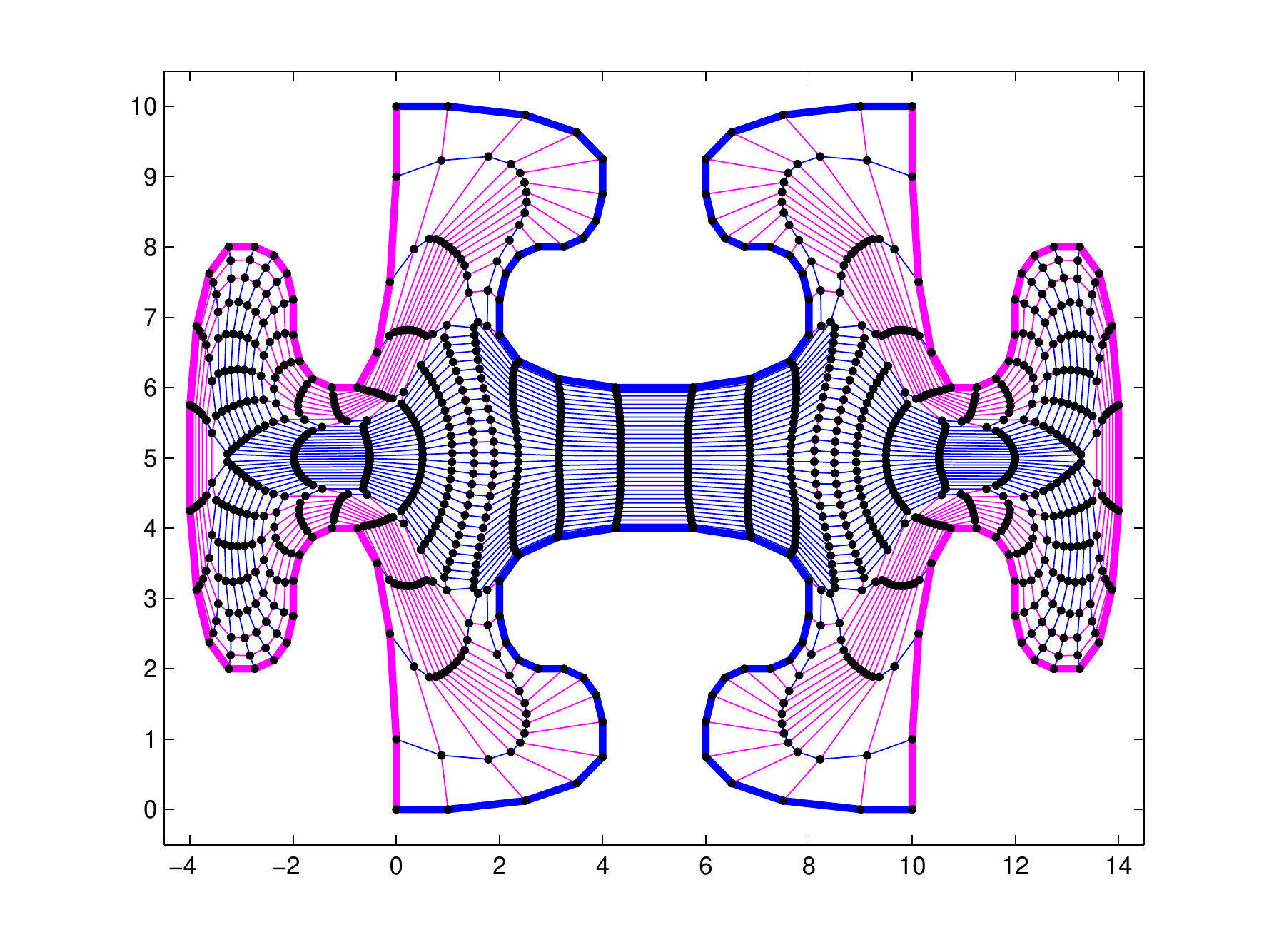}
\caption{Left: Jigsaw puzzle region $\Omega$. Blue curves are the ``south" and ``north" quadratic B-spline boundary curves $\mathbf{F}(\xi,0)$ and  $\mathbf{F}(\xi,1)$ respectively. Pink curves are the ``west" and ``east" quadratic B-spline boundary curves $\mathbf{F}(0,\eta)$ and $\mathbf{F}(1,\eta)$ respectively. Right: the $34 \times 34$ control mesh of the biquadratic B-spline parametrization $\mathbf{F}(\xi,\eta)$ of $\Omega$.}
\label{Fig:Puzzlecont}
\end{figure}

In this section we solve the Poisson equation in the jigsaw puzzle region shown in Figure \ref{Fig:Puzzlecont} (left) and introduced in \cite{Grav14}. We  parametrize this region using an injective biquadratic map $\mathbf{F}(\xi,\eta)$ with a {\it uniform} sequence of knots in both directions and a mesh of $34 \times 34$ control points. These points are computed as the vertices of a quadrilateral mesh \cite{Abe18}, see Figure \ref{Fig:Puzzlecont}, right.

\begin{figure}[htb]
\begin{center}
\includegraphics [scale=0.33]{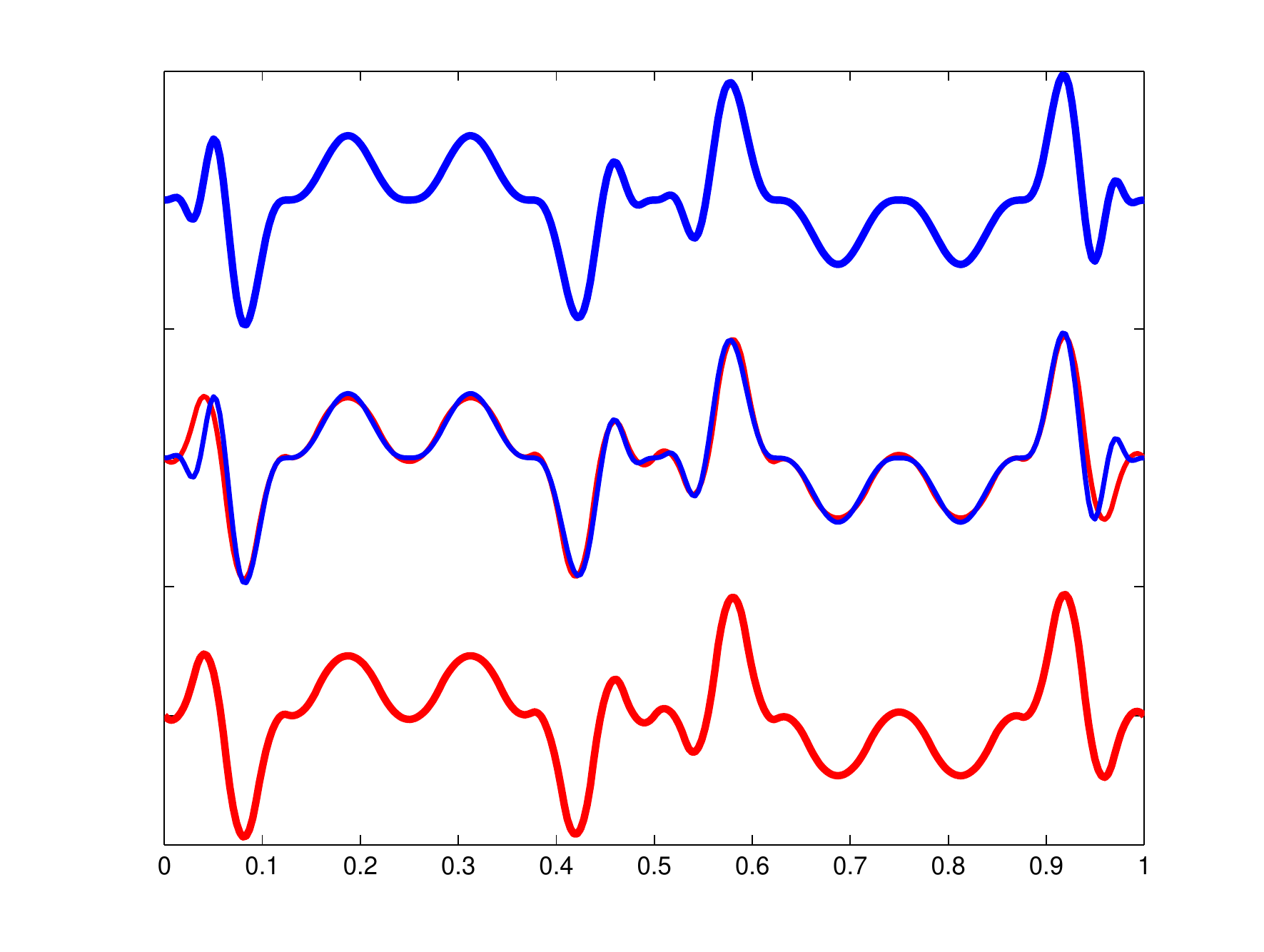}
\includegraphics [scale=0.33]{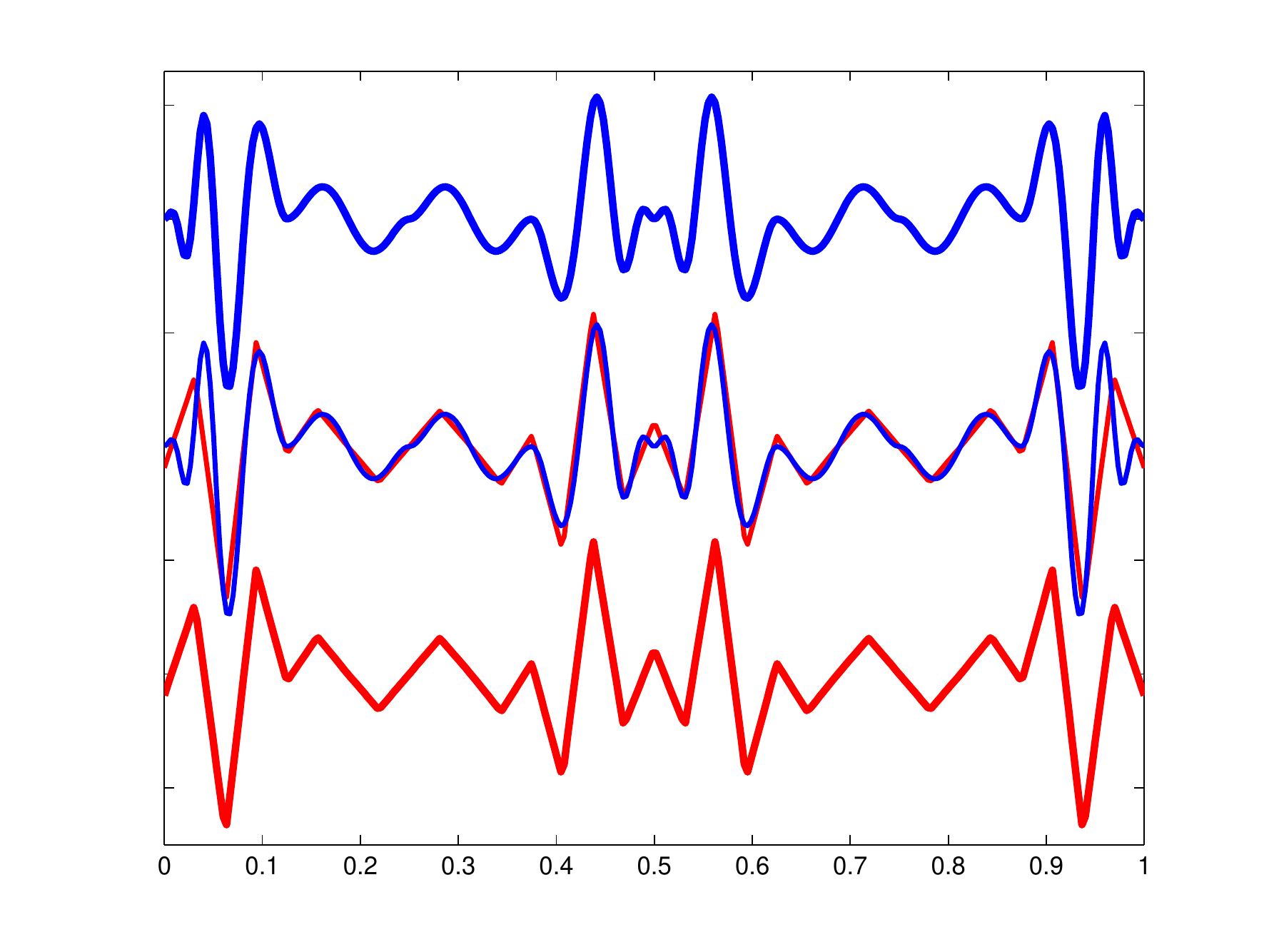}
\end{center}
\caption{ Left top: function $g(\mathbf{F}(\xi,0))$, left bottom: function $u_g^h(\mathbf{F}(\xi,0))$ with uniform knots and $34$ coefficients, left center: comparison of both. Right top: function  $\frac{d}{d \xi}g(\mathbf{F}(\xi,0))$, right bottom: function $\frac{d}{d \xi}u_g^h(\mathbf{F}(\xi,0))$, right center: comparison of both.}
\label{Fig: Curvasg_y_BsplineFronteraSur}
\end{figure}

The function $u_g^h$ is computed interpolating the function $g$ as explained in section \ref{seccionug}. In figure \ref{Fig: Curvasg_y_BsplineFronteraSur} (left) we show the restriction to the ``south" boundary of $\Omega$ of $g$ and its quadratic B-spline approximation $u_g^h$. It is clear that the approximation is good, except in the middle and the extremes. This is better observed in Figure \ref{Fig: Curvasg_y_BsplineFronteraSur} right, where we compare the derivative of $g$ with the derivative of $u_g^h$ ( a piecewise linear function) both restricted to the ``south" boundary of $\Omega$. In this figure we observe that the derivative of the quadratic B-spline function is not able to represent faithfully the frequencies and amplitudes of the derivative of $g$. In other words, to obtain a better approximation of the derivative of $g$ we need a B-spline quadratic function with more degrees of freedom.

To get extra degrees of freedom new knots must be inserted in the sequence $t^{\xi}$ in the subintervals where the error $\frac{d}{d \xi}g(\mathbf{F}(\xi,0))-\frac{d}{d \xi} u_g^h(\mathbf{F}(\xi,0))$ is big. More precisely we insert 6 knots in the interval $(0,0.1)$ and $(0.9,1)$ and 7 knots in $(0.4,0.6)$. The corresponding B-spline quadratic function has now $53$ degrees of freedom. The same procedure is repeated for the ``north" boundary of
$\Omega$.
%

\begin{figure}[htb]
\begin{center}
\includegraphics [scale=0.35]{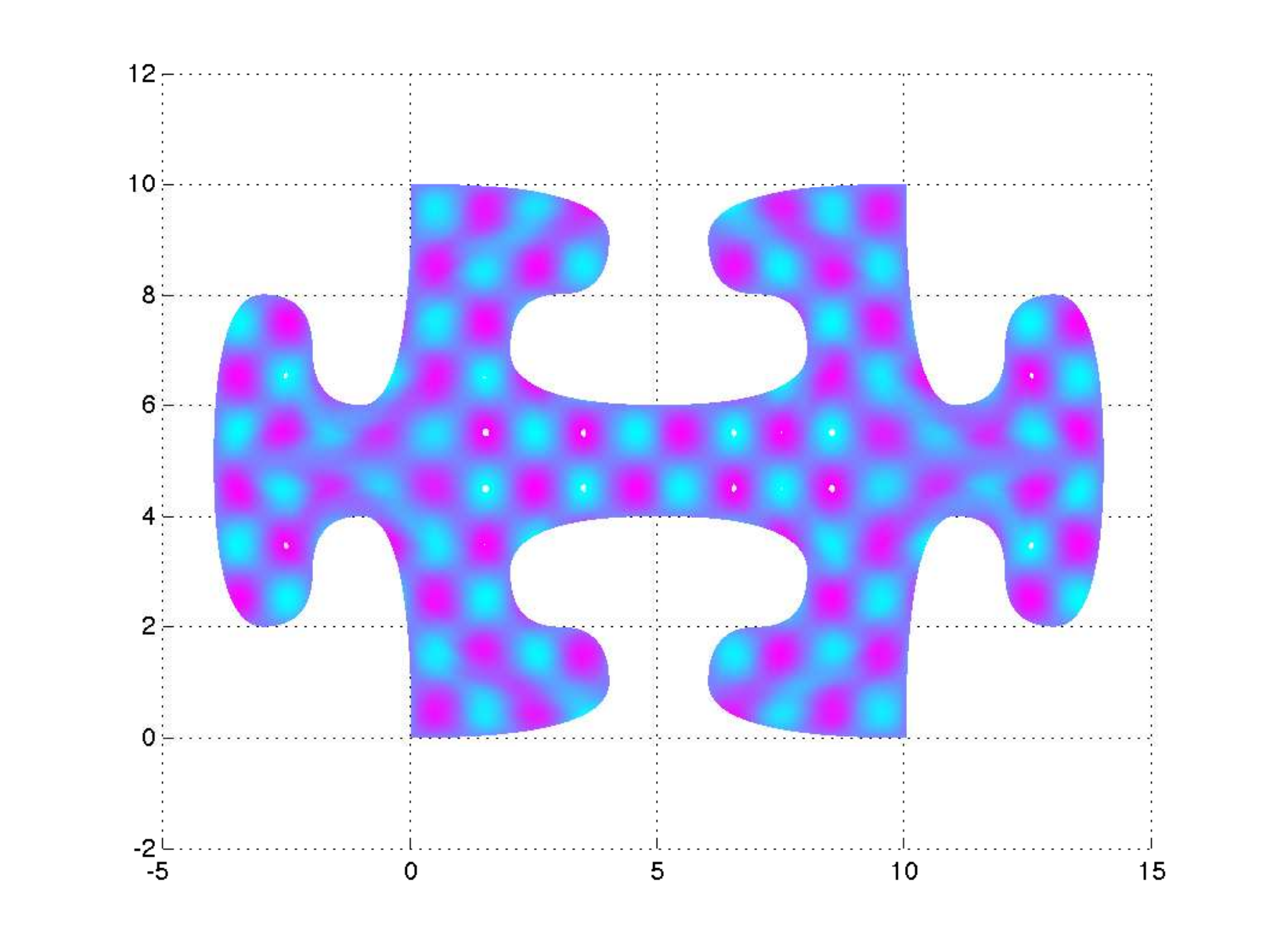}
\includegraphics [scale=0.35]{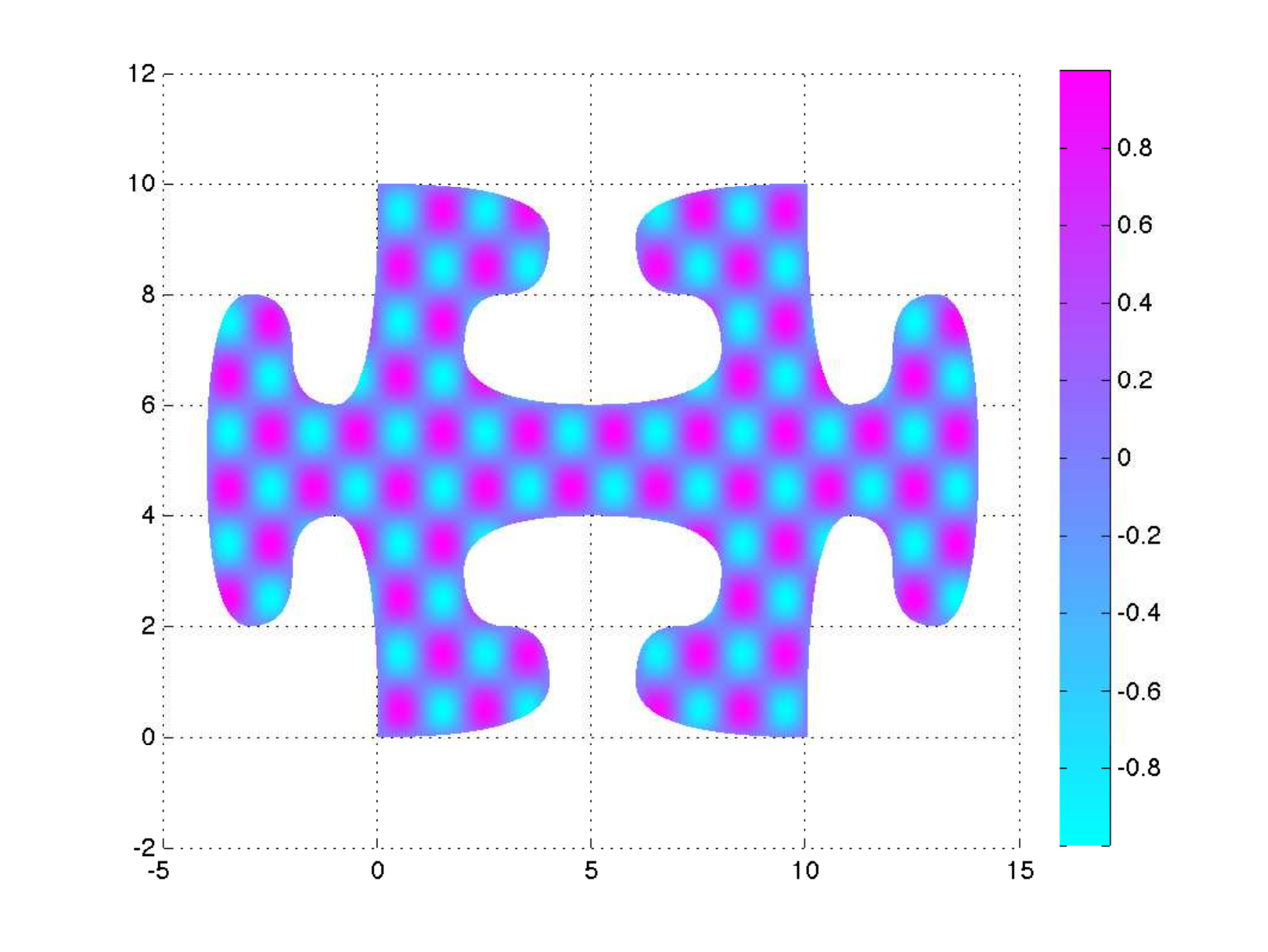}
\end{center}
\caption{Left: 2D view of the approximated solution of Poisson equation using biquadratic B-splines with $53 \times 53$ control points. Right: 2D view of the exact solution.}
\label{Fig:Vantroi53}
\end{figure}

The sequence of knots $t^{\xi}$ obtained after inserting the new knots is also used as $t^{\eta}$. Therefore, the new space $\mathbb{S}_{3,t^\xi} \bigotimes \mathbb{S}_{3,t^\eta}$ has dimension $53 \times 53$. Even when the map $\mathbf{F}(\xi,\eta)$ is the same, its control points in the new basis must be computed \cite{Abe18}. The approximated solution $u^h$ is a biquadratic B-spline function with $53 \times 53$ degrees of freedom, which are computed solving the linear system of section 3.  In Figure \ref{Fig:Vantroi53} we show a 2D view of the approximated solution $u^h$ and the exact solution $u$.

\begin{figure}[htb]
\begin{center}
\includegraphics [scale=0.3]{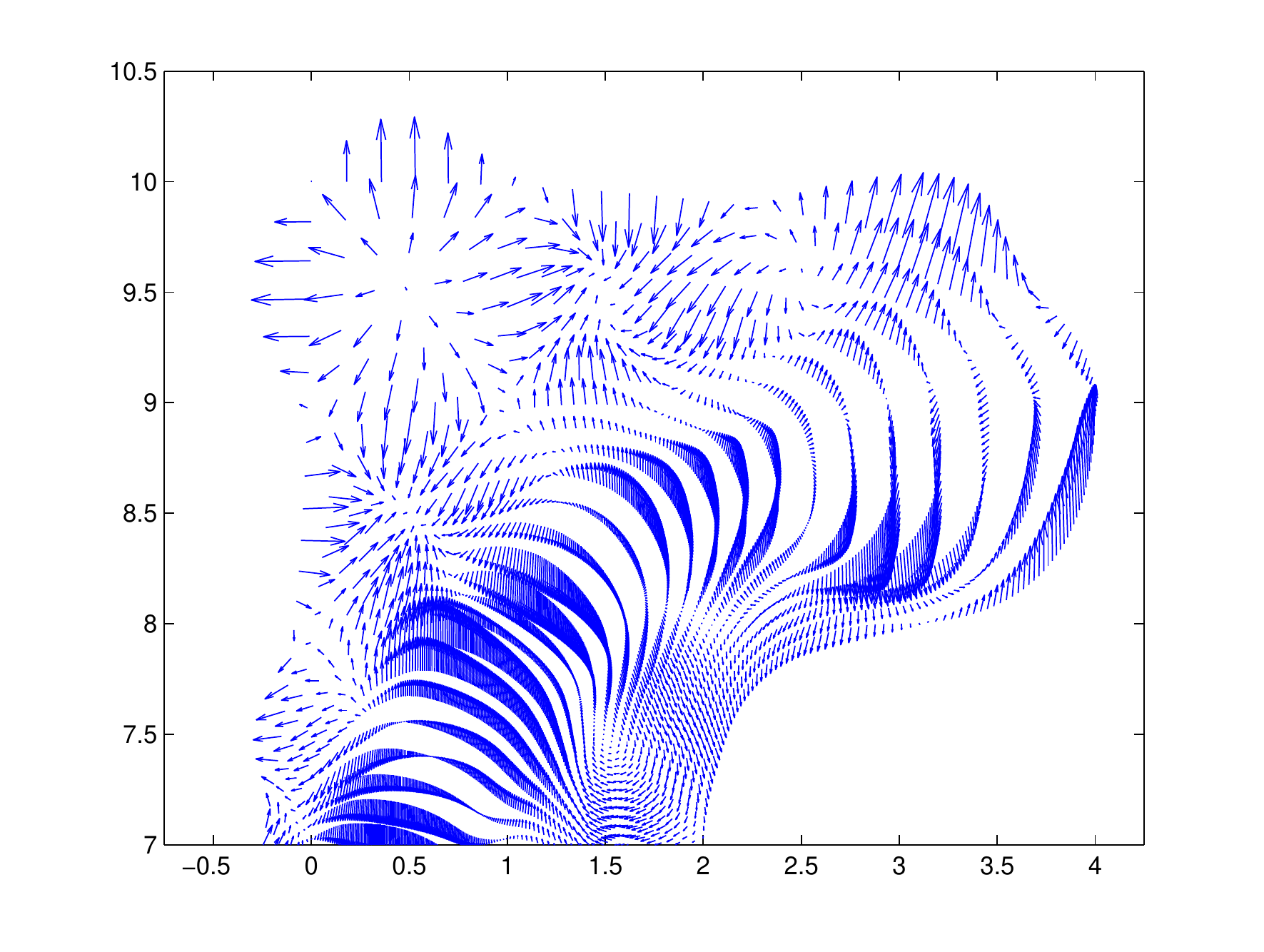}
\includegraphics [scale=0.3]{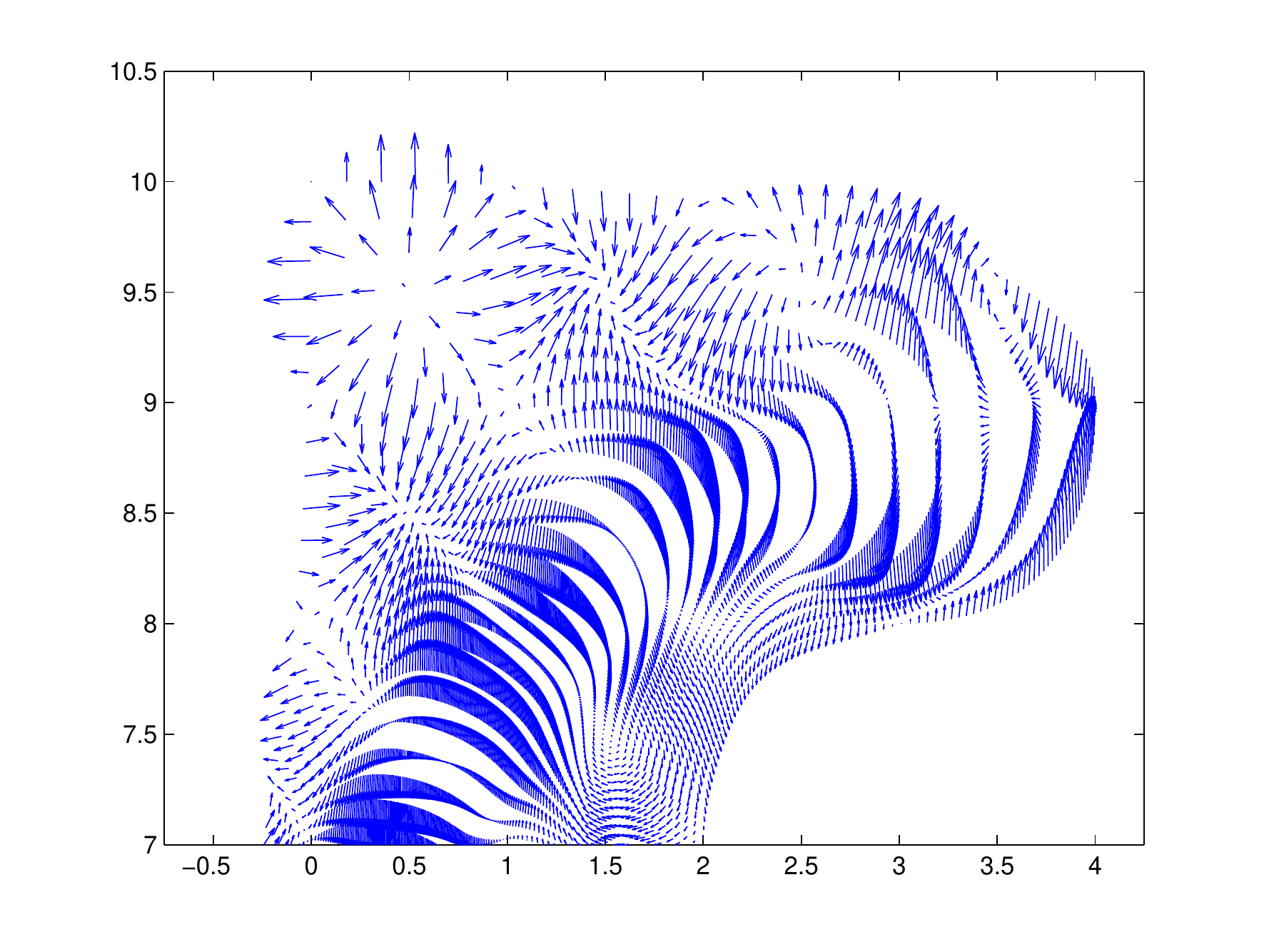}
\end{center}
\caption{Zoom of the vector field. Left: vector field for the quadratic B-spline approximation with $53 \times 53$ control points. Right: vector field of the exact solution.}
\label{Fig:Vantroi53CampVec}
\end{figure}

The error $L_2$ error of $u^h$ given by (\ref{l2error}) is equal to $0.1462$, a  relative small value, but the $H_1$ error is $18.5486$. This can be better observed in Figure \ref{Fig:Vantroi53CampVec}, where we show a zoom of the exact and the approximated vector field in a section of the physical domain $\Omega$.

In order to obtain a better approximation of the vector field, we increase the dimension of the approximating spline space, inserting a new knot in the middle point between two consecutive knots. The new space has dimension $104 \times 104$. The map $\mathbf{F}(\xi,\eta)$ is still the same, but its $104 \times 104$ control points in the new basis must be computed. Moreover, we compute the approximated solution $u^h$, which also belongs to the same biquadratic B-spline space. The result is a better approximation: the $L_2$ error is now $0.0126$ while the $H_1$ error is $3.6864$.

\subsection{Poisson equation with discontinuous gradient}

In this section we solve the Poisson equation
$$ -\triangle u(x,y) = f(x,y), \;\;(x,y) \in \Omega$$
with Dirichlet boundary condition (\ref{Hebc}). The function $f(x,y)$ is computed in such away that the exact solution $u(x,y)$ is given by,
\begin{equation}
\label{funMarina01}
\begin{array}{ccc}
u(x,y)=\exp\left( \alpha   \sqrt{(x-x_0)^2 + (y-y_0)^2}\right)
          &+&\exp\left( \beta     \sqrt{(x-x_1)^2 + (y-y_1)^2}\right)+\\
          &+&\exp\left( \gamma\sqrt{(x-x_2)^2 + (y-y_2)^2}\right)
\end{array}
\end{equation}

where the real values $\alpha,\beta$ and $\gamma$ and the points $(x_0,y_0),\ (x_1,y_1)$ and $(x_2,y_2)$ are known. The function $g(x,y)$ is the restriction of $u(x,y)$ to the boundary of $\Omega$. This problem is solved in \cite{Bro16} on the unit square $[0,1]^2$. Here we solve it on several irregular regions. The main difficulty is the discontinuity of the gradient of $u(x,y)$ in the points $(x_0,y_0),\ (x_1,y_1)$ and $(x_2,y_2)$.
%

In the following experiments $\alpha=\beta=\gamma=7$ and the points involved in the description of the problem and in its solution (\ref{funMarina01}) are: $(x_0,y_0)=\mathbf{F}(\xi_a,\eta_a)$, $(x_1,y_1)=\mathbf{F}(\xi_b,\eta_b)$ and $(x_2,y_2)=\mathbf{F}(\xi_c,\eta_c)$, where $\xi_a=\eta_a=0.25$, $\xi_b=\eta_b=0.5$ and  $\xi_c=\eta_c=0.75$. The first step to obtain the approximated solution $u^h$ is to compute the B-spline biquadratic parametrization $\mathbf{F}(\xi,\eta)$ of the physical region $\Omega$. The sequences of knots $t^{\xi}$ and $t^{\eta}$ that we use to define the space of biquadratic splines are {\it nonuniform}. More precisely, the distribution of knots in $t^{\xi}$ is more concentrated near the parametric values $\xi_a,\xi_b$ and  $\xi_c$. Similarly, the sequence $t^{\eta}$ contains more knots near the parametric values  $\eta_a,\eta_b$ and  $\eta_c$.

In table \ref{tablaMarinaSuave} we show the results for different regions. The number of degrees of freedom used to compute $u^h$ is reported in the second column of the table. The other two columns contain the errors (\ref{l2error}) and (\ref{H1error}). As we observe, the $L_2$ error oscillates between $10^{-2}$ and $10^{-4}$, but the $H_1$ error is approximately two orders bigger. It means that $u^h$ could be considered as a good approximation of the exact solution $u$, but partial derivatives of $u^h$ are not good approximations of partial derivatives of $u$.

\begin{table}[hbt]
\begin{center}
\begin{tabular}{||c|c|c|c||}
\hline
\hline
Region                & Degrees of freedom & $L_2\;error$ & $H_1\;error$\\\hline\hline
Havana bay            &  $116 \times 110$         &  0.0648      &   6.8375\\\hline
Toba lake             &  $172 \times 172$         &  8.5626e-4   &   0.0577\\\hline
Gibraltar channel     &  $ 96 \times 112$         & 0.0190       &   2.4041\\\hline
Grijalva channel      &  $124 \times 44 $         & 0.0055       &   2.3928\\\hline
P\'atzcuaro lake      &  $108 \times 108$         & 9.0471e-4    &   0.2075\\\hline
V. de Bravo reservoir &  $156 \times 156$         & 4.5784e-4    &   0.0328\\
\hline
\hline
\end{tabular}
\caption{\label{tablaMarinaSuave} Errors of the biquadratic B-spline solution of Poisson equation with exact solution (\ref{funMarina01}) on several physical regions. The parameters $\xi_a,\:\xi_b,\:\xi_c$ are simples knots in $t^{\xi}$ and the parameters $\eta_a,\:\eta_b,\:\eta_c$ are simple knots in $t^{\eta}$.}
\end{center}
\end{table}

\begin{figure}[hbt]
\begin{center}
\includegraphics [scale=0.3]{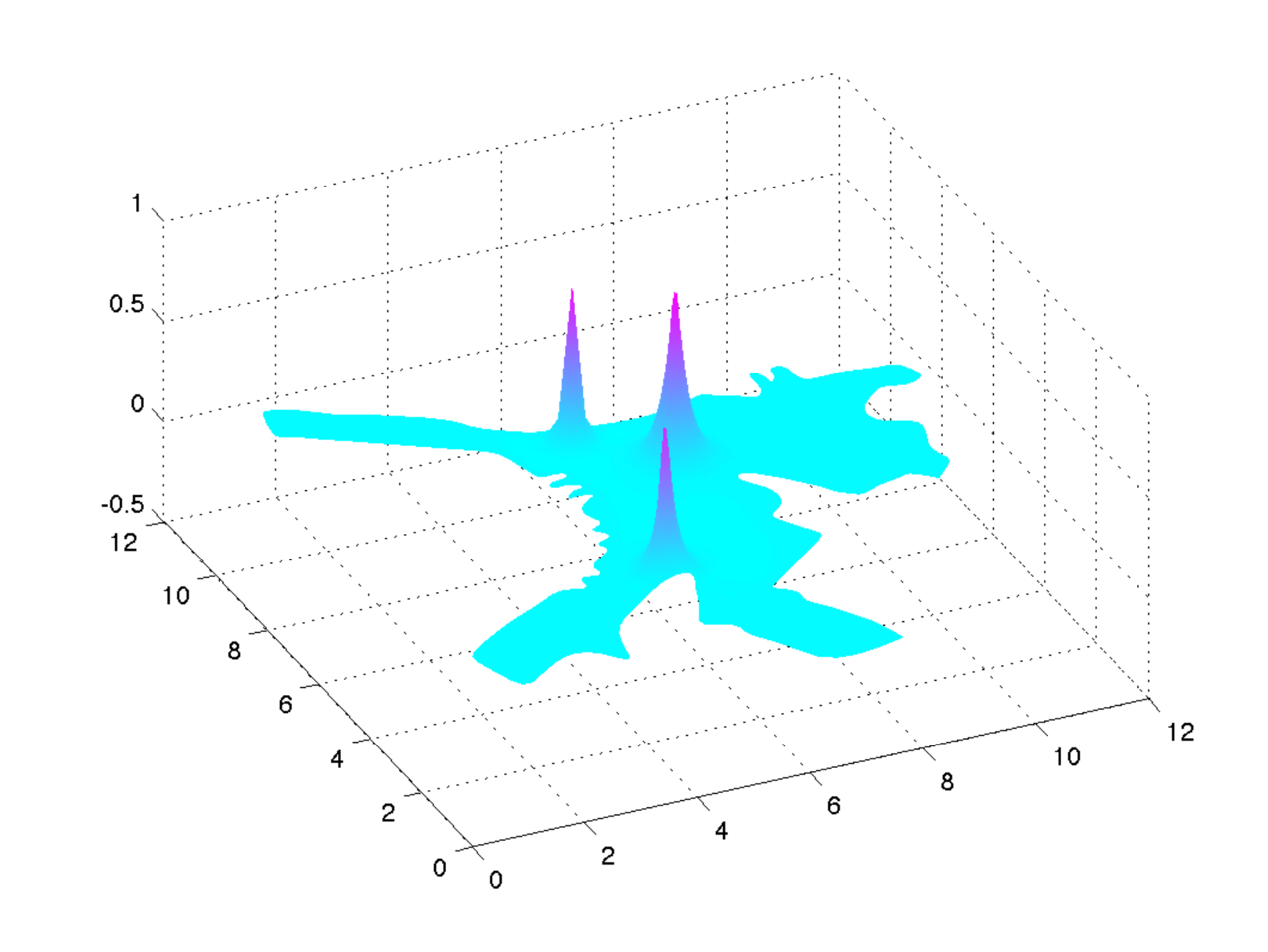}
\includegraphics [scale=0.3]{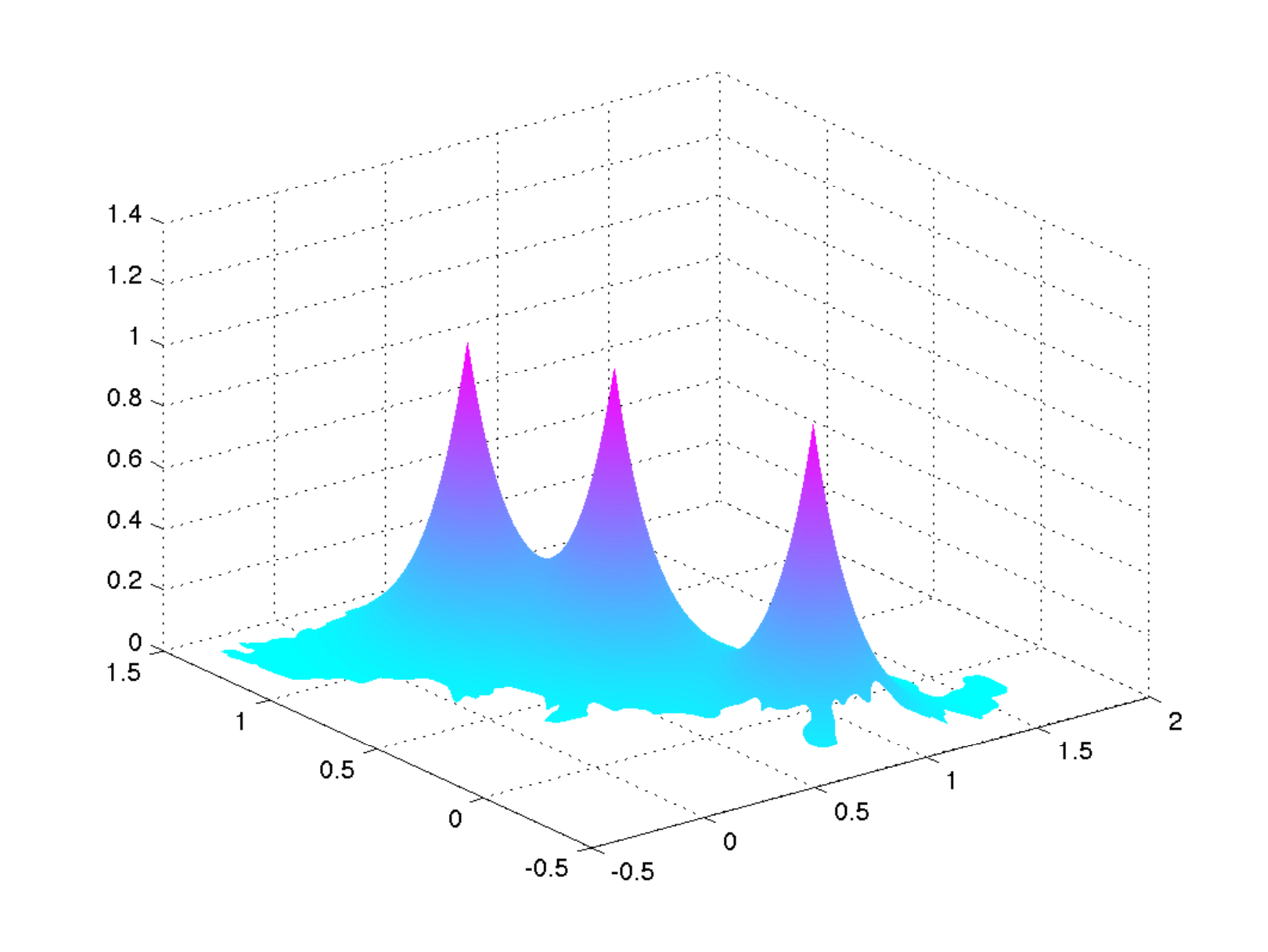}\\
\includegraphics [scale=0.3]{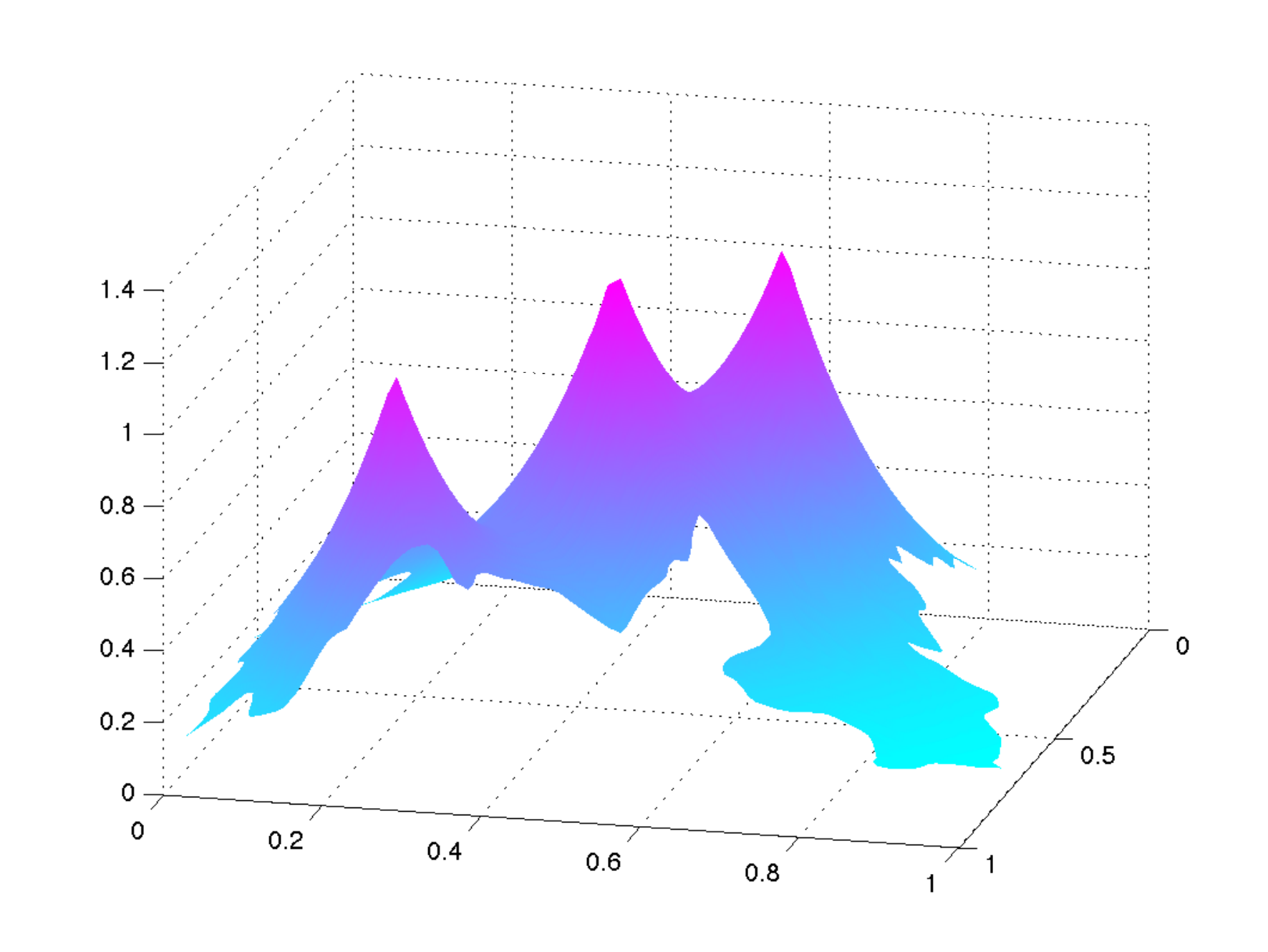}
\includegraphics [scale=0.3]{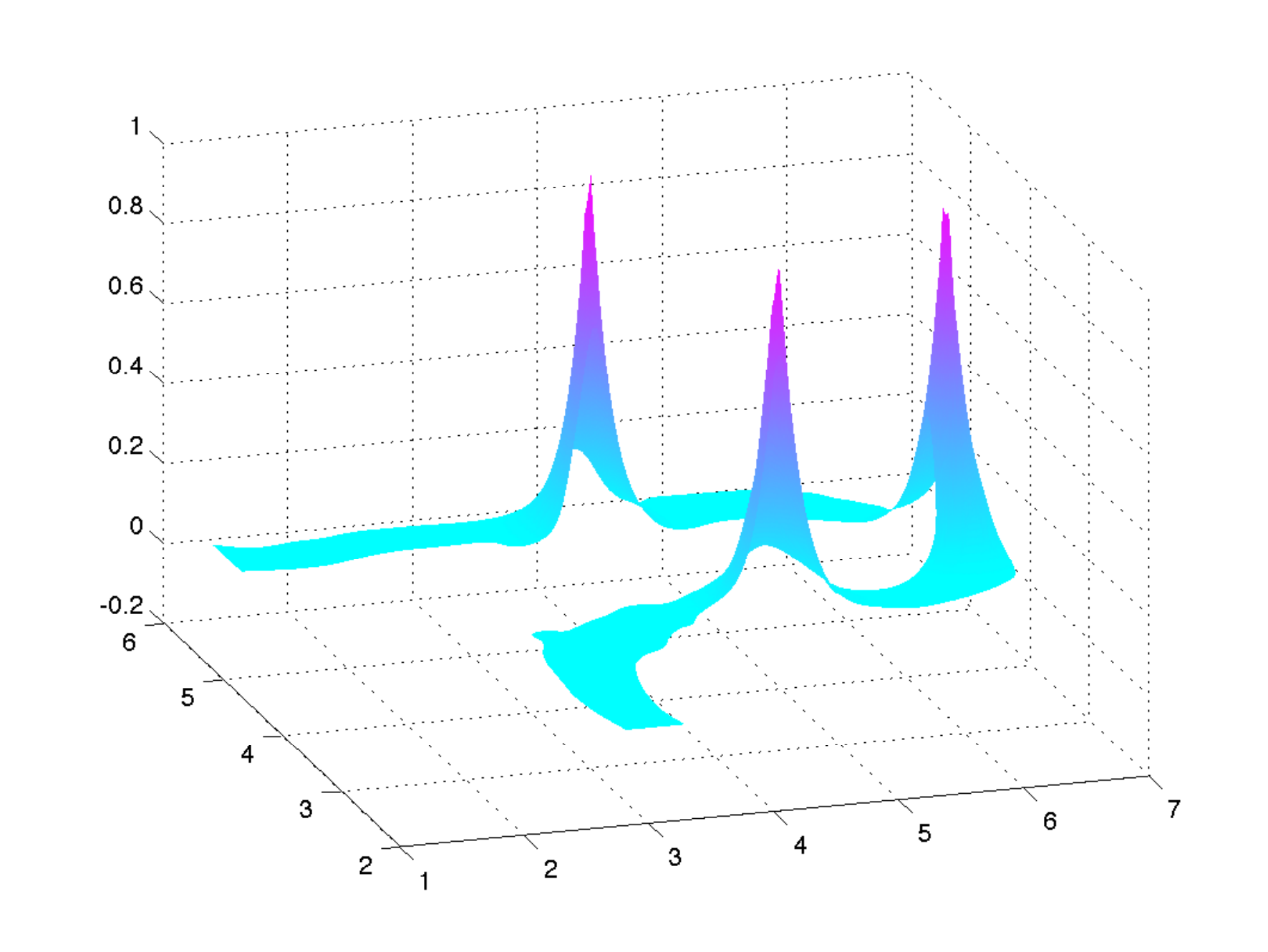}
\end{center}
\caption{ Biquadratic B-spline functions $u^h$ approximating the exact solution (\ref{funMarina01}) of Poisson equation for several regions. Top left: Havana bay, top right: Toba lake. Bottom left: V. de Bravo reservoir, bottom right: Grijalva channel.}
\label{Fig:SolMarinaEjBSp}
\end{figure}

In Figure \ref{Fig:SolMarinaEjBSp} we show the biquadratic B-spline functions $u^h$ for some of the physical regions reported in table \ref{tablaMarinaSuave}. We recall that the B-spline basis functions $B_i^3(\xi)$ and $B_j^3(\eta)$ used to construct the approximated solution $u^h$ are $C^1$ continuous, since the corresponding sequences of knots $t^{\xi}$ and $t^{\eta}$ are composed by simple knots. However, the gradient of the exact solution (\ref{funMarina01}) is not defined in three points. Hence, the smooth B-spline solution $u^h$ approximates the exact solution $u$ in these points, but partial derivatives of $u^h$ are not good approximations of partial derivatives of $u$. To overcome this difficulty, we include two times the parametric values $\xi_a,\:\xi_b$ and $\xi_c$ in the sequence of knots  $t^{\xi}$, and similarly the parametric value $\eta_a,\:\eta_b$ and $\eta_c$ are included two times in the sequence $t^{\eta}$. The corresponding B-spline functions are only $C^0$ continuous in these points, therefore $u^h$ approximates better the exact solution $u$. Observe that the introduction of new knots does not change the map $\mathbf{F}(\xi,\eta)$, which is still differentiable, but the expression of $\mathbf{F}(\xi,\eta)$ in the new basis must be computed.

\begin{figure}[hbt]
\begin{center}
\includegraphics [scale=0.26]{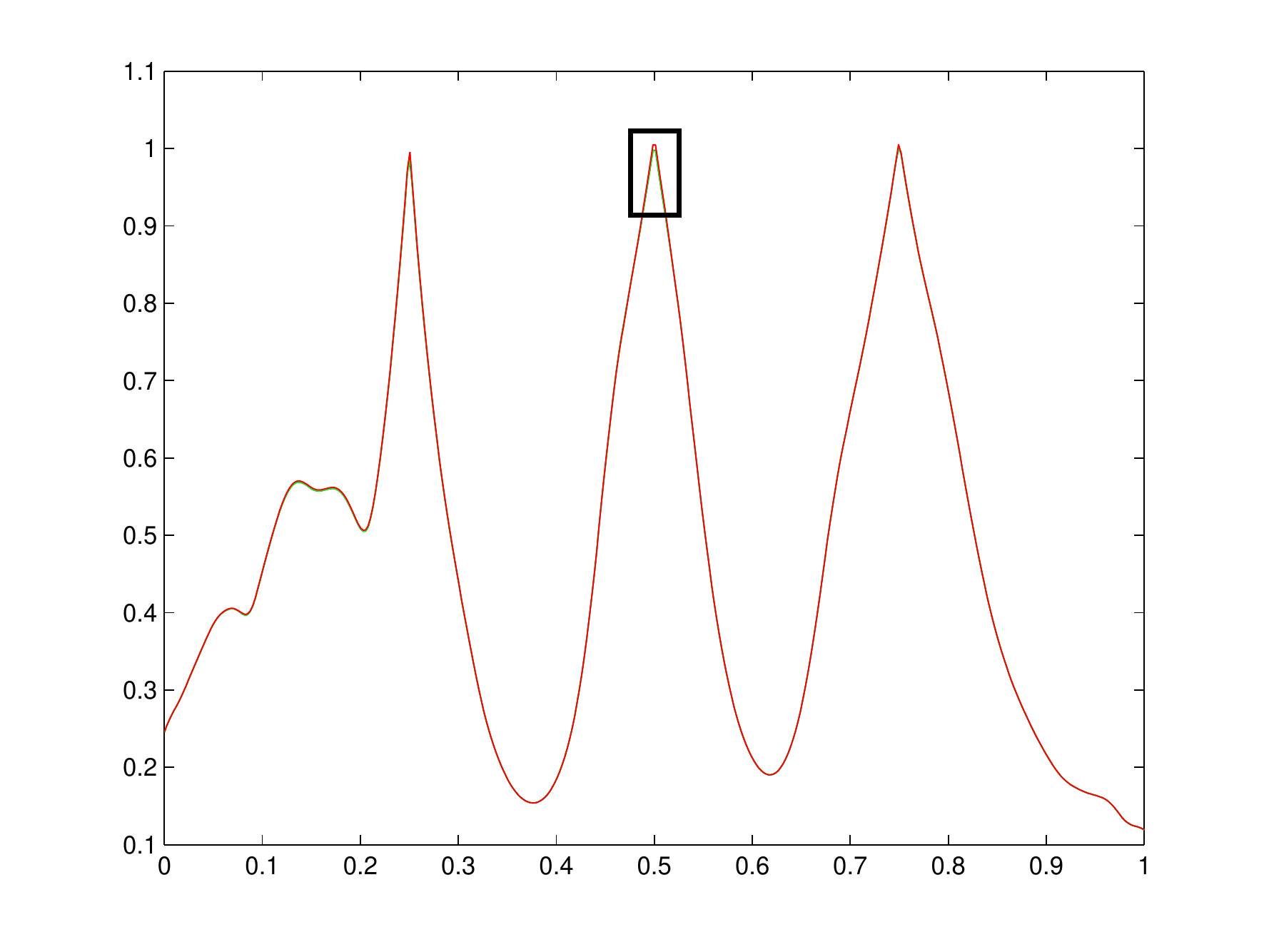}
\includegraphics [scale=0.26]{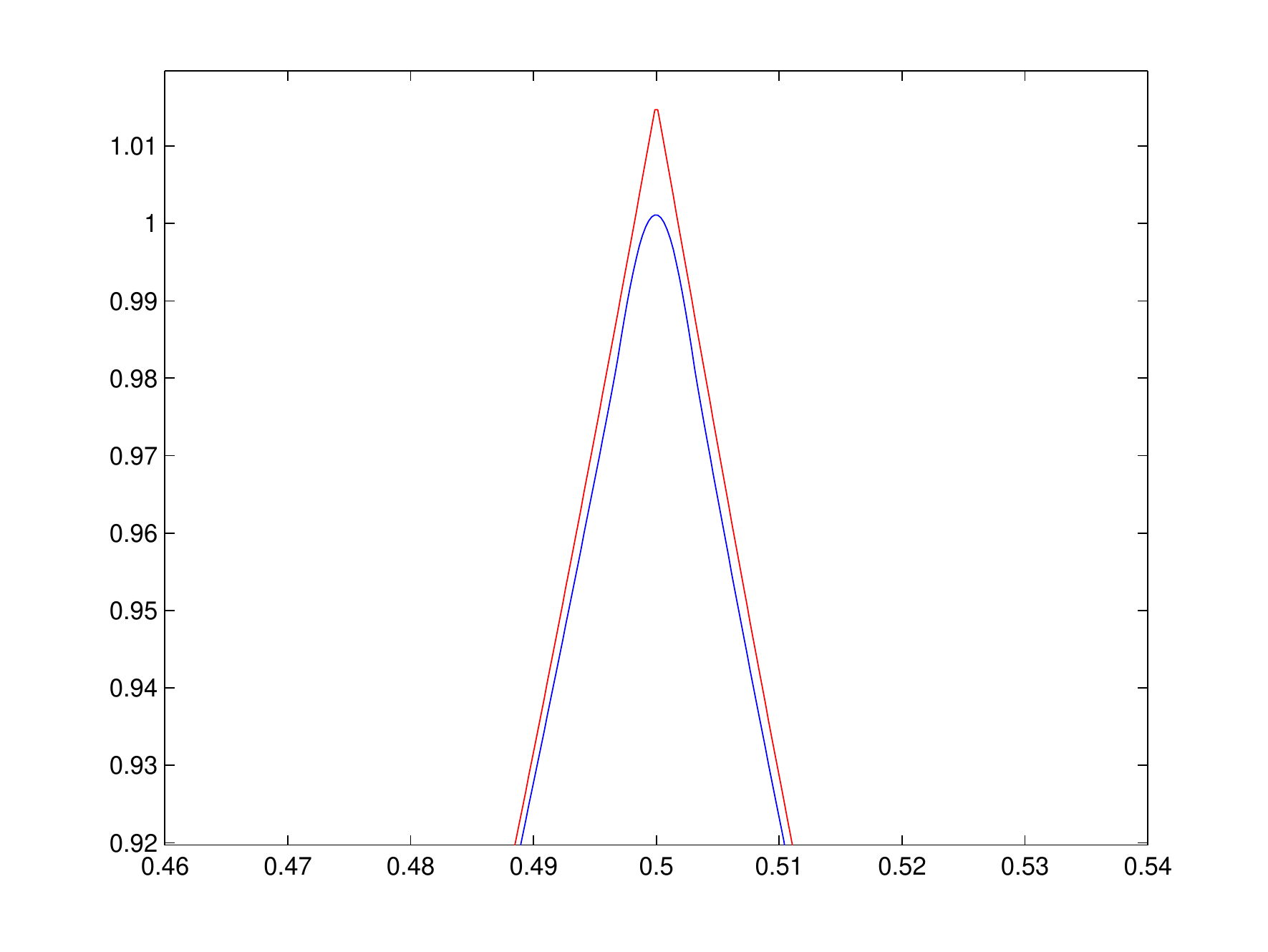}
\includegraphics [scale=0.26]{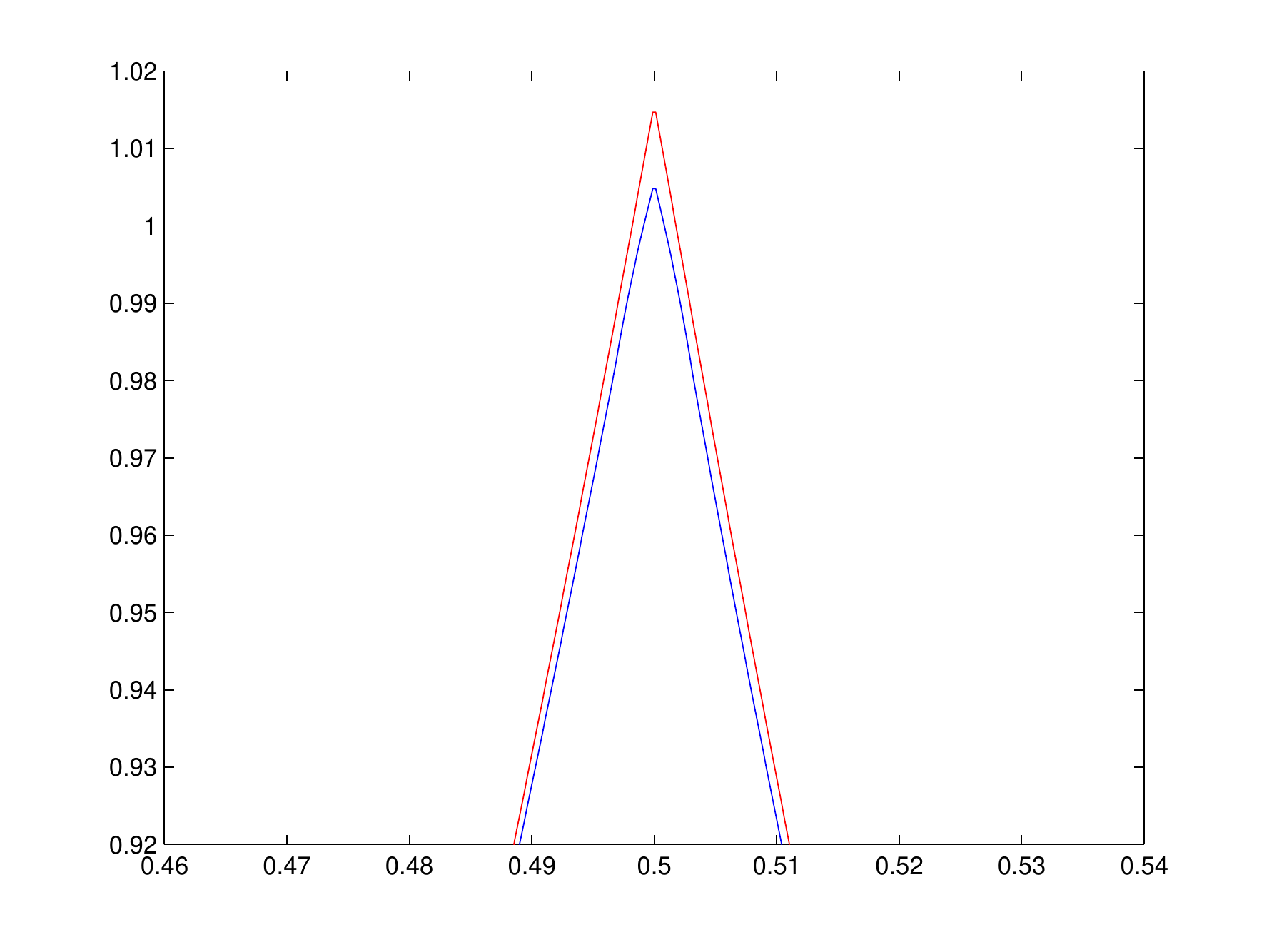}
\end{center}
\caption{Left: graph of $u(\mathbf{F}(\xi,\xi))$ with $u$ given by (\ref{funMarina01}). Center and right: zoom of $u^h(\mathbf{F}(\xi,\xi))$ (in blue) on the P\'atzcuaro lake in comparison with $u(\mathbf{F}(\xi,\xi))$ ( in red). In center graph $u^h$ is computed with simple knots, in right graph $u^h$ has repeated knots.}
\label{Fig:PatzSolsMarina1D}
\end{figure}
In Figure \ref{Fig:PatzSolsMarina1D} we compare the graph of the functions $u(\mathbf{F}(\xi,\xi))$ and $u^h(\mathbf{F}(\xi,\xi))$, where $\mathbf{F}$ is the parametrization of P\'atzcuaro lake and $u^h$ is the biquadratic B-spline approximation to the exact solution $u$. Observe that these curves contain the singular points of $u$. The left image shows the function  $u(\mathbf{F}(\xi,\xi)),\;0 \leq \xi \leq 1$, while center and right images show a zoom of $u(\mathbf{F}(\xi,\xi))$ and $u^h(\mathbf{F}(\xi,\xi))$ restricted to the black rectangle in the left image. This rectangle contains the point $\xi_b=0.5$.  The center graph shows in blue the approximated B-spline solution $u^h$ obtained for a sequence of simple knots. We observe that the exact solution $u(\mathbf{F}(\xi,\xi))$ (in red) is not differentiable in $\xi_b$, but $u^h(\mathbf{F}(\xi,\xi))$ has continuous derivative in this point. The right graph shows in blue the approximated B-spline solution $u^h$ corresponding to a sequence of knots  $t^{\xi}$, where $\xi_b=0.5$ is repeated,  and a sequence of knots $t^{\eta}$,  where $\eta_b=0.5$ is also repeated. The result is that $u^h$ has the same behavior that $u$ since it is not differentiable in $\xi_b=0.5$. The effect of repeating the knots $\xi_a,\xi_b$ and $\xi_c$ in $t^{\xi}$ and $\eta_a,\eta_b$ and $\eta_c$ in $t^{\eta}$ is shown in table \ref{tablaMarinaNoSuave}, which contains the results for the same physical regions of table \ref{tablaMarinaSuave}.

\begin{table}[hbt]
\begin{center}
\begin{tabular}{||c|c|c|c||}
\hline
\hline
Region                & Degrees of freedom & $L_2\;error$  & $H_1\;error$\\\hline\hline
Havana bay            &  $119 \times 113$         & 0.0084        &   1.0420\\\hline
Toba  lake            &  $175 \times 175$         & 8.5733e-4     &   0.0188\\\hline
Gibraltar channel     &  $ 99 \times 115$         & 0.0190        &   1.2613\\\hline
Grijalva channel      &  $127 \times 47 $         & 0.0054        &   1.7532\\\hline
P\'atzcuaro lake      &  $111 \times 111$         & 9.0369e-4     &   0.1973\\\hline
V. de Bravo reservoir &  $159 \times 159$         & 4.5737e-4     &   0.0187\\
\hline
\hline
\end{tabular}
\caption{\label{tablaMarinaNoSuave} Errors of the biquadratic B-spline solution of Poisson equation with exact solution (\ref{funMarina01}). The parameters $\xi_a,\:\xi_b,\:\xi_c$ are double knots in $t^{\xi}$ and $\eta_a,\:\eta_b,\:\eta_c$ are double knots in $t^{\eta}$.}
\end{center}
\end{table}

\begin{figure}[hbt]
\begin{center}
\includegraphics [scale=0.22]{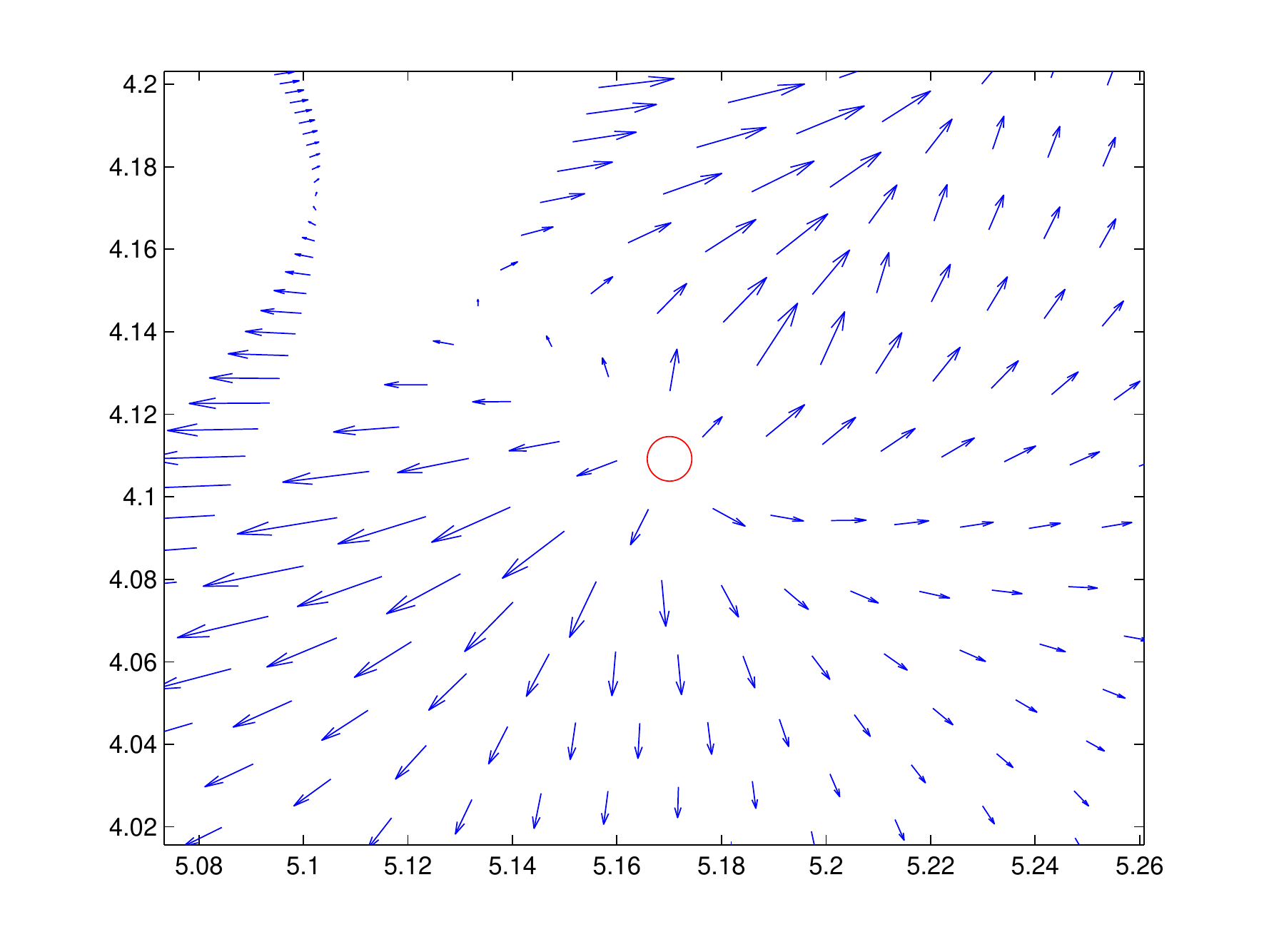}
\includegraphics [scale=0.22]{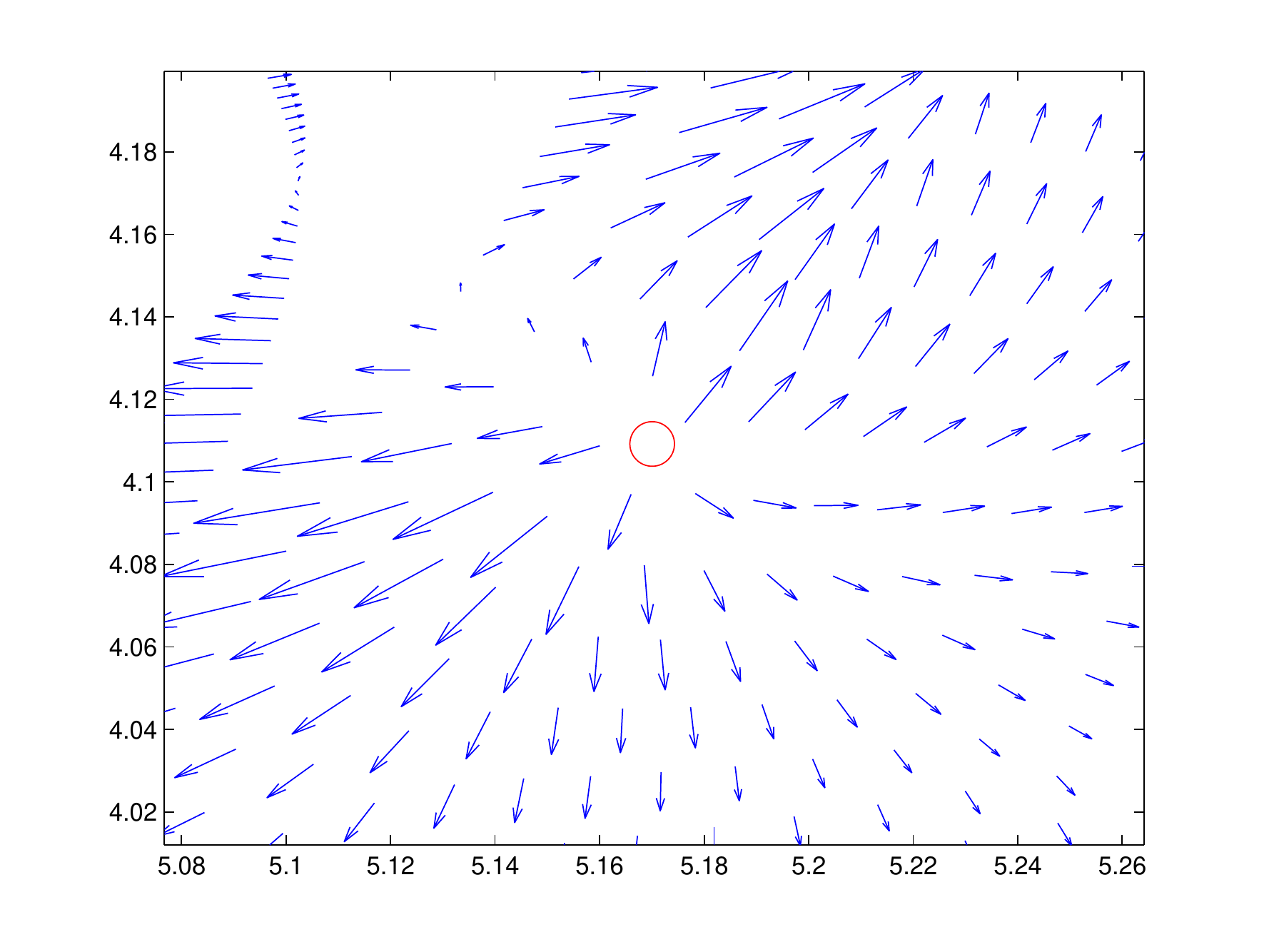}
\includegraphics [scale=0.22]{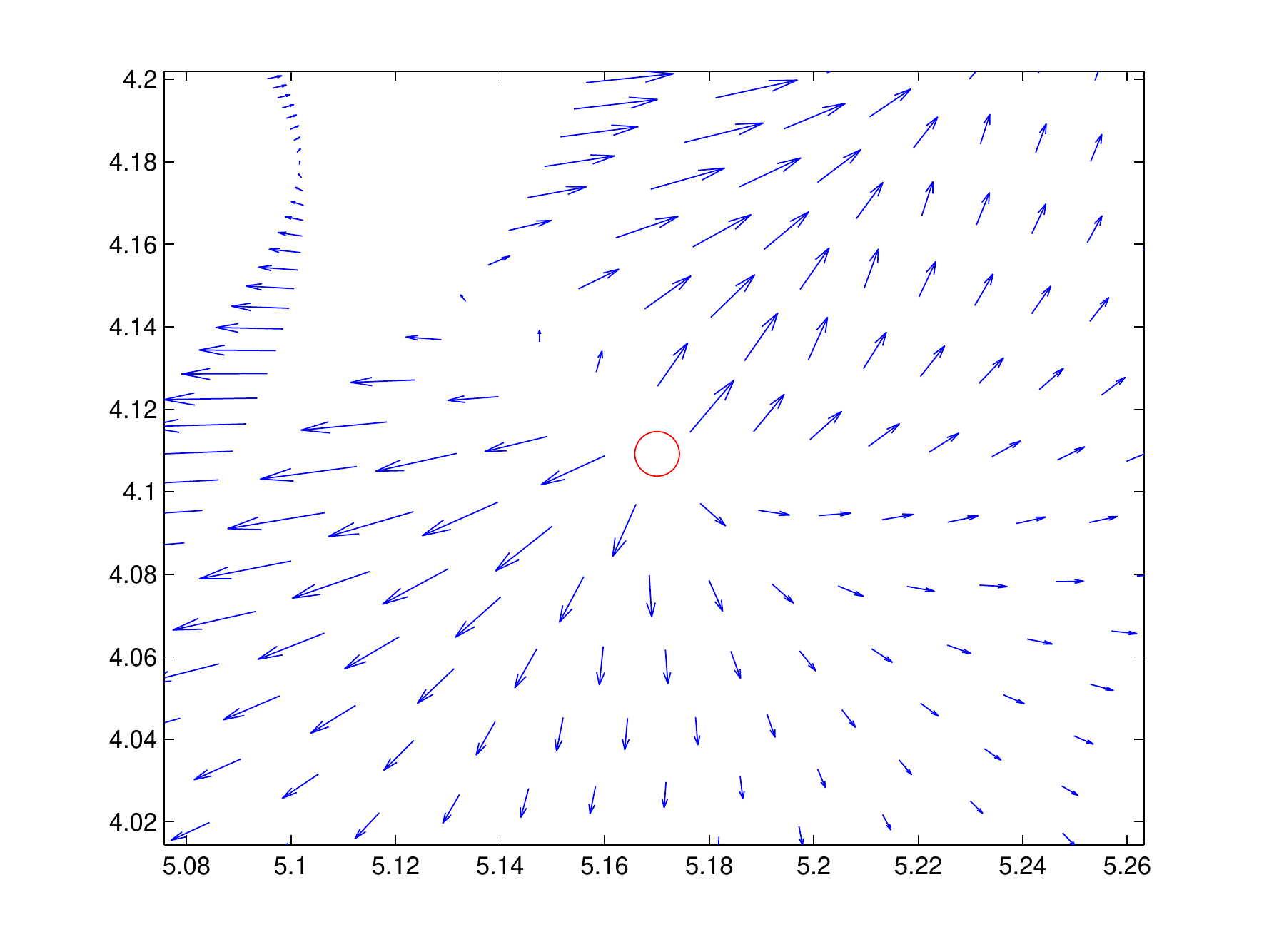}
\end{center}
\caption{Zoom of the vector field of the exact solution and the biquadratic B-spline solution of Poisson equation. Left: for the smooth B-spline approximation with simple knots. Center: for the B-spline approximation with double knots. Right: for the exact solution.}
\label{Fig:CampVecMarinaBahHab}
\end{figure}

Comparing tables \ref{tablaMarinaSuave}  and \ref{tablaMarinaNoSuave} we observe that in each parametric direction, the number of degrees of freedom is increased in 3, because we repeat 3 knots in the corresponding sequences $t^{\xi}$ and $t^{\eta}$. As a consequence, the $L_2$ and the $H_1$ errors are reduced in general. The reduction is significative for the $H_1$ error, since repeating knots we obtain a better approximation of the vector field of the exact solution. This is illustrated in Figure \ref{Fig:CampVecMarinaBahHab}, where we show the vector field near a singular point for Havana bay. The left and center images of this figure show the vector field of the biquadratic B-spline function $u^h$ for simple and repeated knots respectively. The right image shows the vector field of the exact solution $u$. It is easy to see that the size of the arrows near the singular point is smaller for the left image, which means that the field is smoother in this point. Moreover, the vector fields for the center and right images are very similar.

\subsection{Helmholtz equation with variable frequency}

The wave function $u(x,y)$ that satisfies a Schr\"{o}dinger equation model of two interacting atoms \cite{Mit13} is
the solution of the Helmholtz equation (\ref{Heq}) with
\begin{equation}
k(x,y)=\frac{1}{(\alpha + r(x,y))}
\label{kSho}
\end{equation}
where $\alpha$ is a parameter $r(x,y)=\sqrt{(x-x_0)^2+(y-y_0)^2}$ and
\begin{equation}
f(x,y)=\frac{(\alpha-r(x,y)) \cos(k(x,y))}{(\alpha+r(x,y))^3r(x,y)} \label{fxy}
\end{equation}
In this case, the exact solution of Helmholtz equation is given by
\begin{equation}
u(x,y)=\sin(k(x,y))
\label{SolHelm}
\end{equation}
The function (\ref{SolHelm}) has discontinuous gradient at $(x_0,y_0)$ and it is highly oscillatory near that point. The number of oscillations $M$ is determined by the parameter $\alpha=\frac{1}{M \pi}$.

%

\subsubsection{Experiments for exact solution with only one oscillation. }
In this section we solve the Helmholtz equation with $k(x,y)$ given by (\ref{kSho}) for several regions with irregular boundary. In all the examples we select  $\alpha=\frac{1}{\pi}$ and we compute the point $(x_0,y_0)$ as $\mathbf{F}(\widetilde{\xi},\widetilde{\eta})$, where $(\widetilde{\xi},\widetilde{\eta})=(0.5,0.5)$ and the parametrization $\mathbf{F}(\xi,\eta)$ is the biquadratic B-spline function given by (\ref{MapBiCuadExpli}), with control points computed as the vertices of a quadrilateral mesh \cite{Abe18}.

\begin{figure}[hbt]
\begin{center}
\includegraphics [scale=0.27]{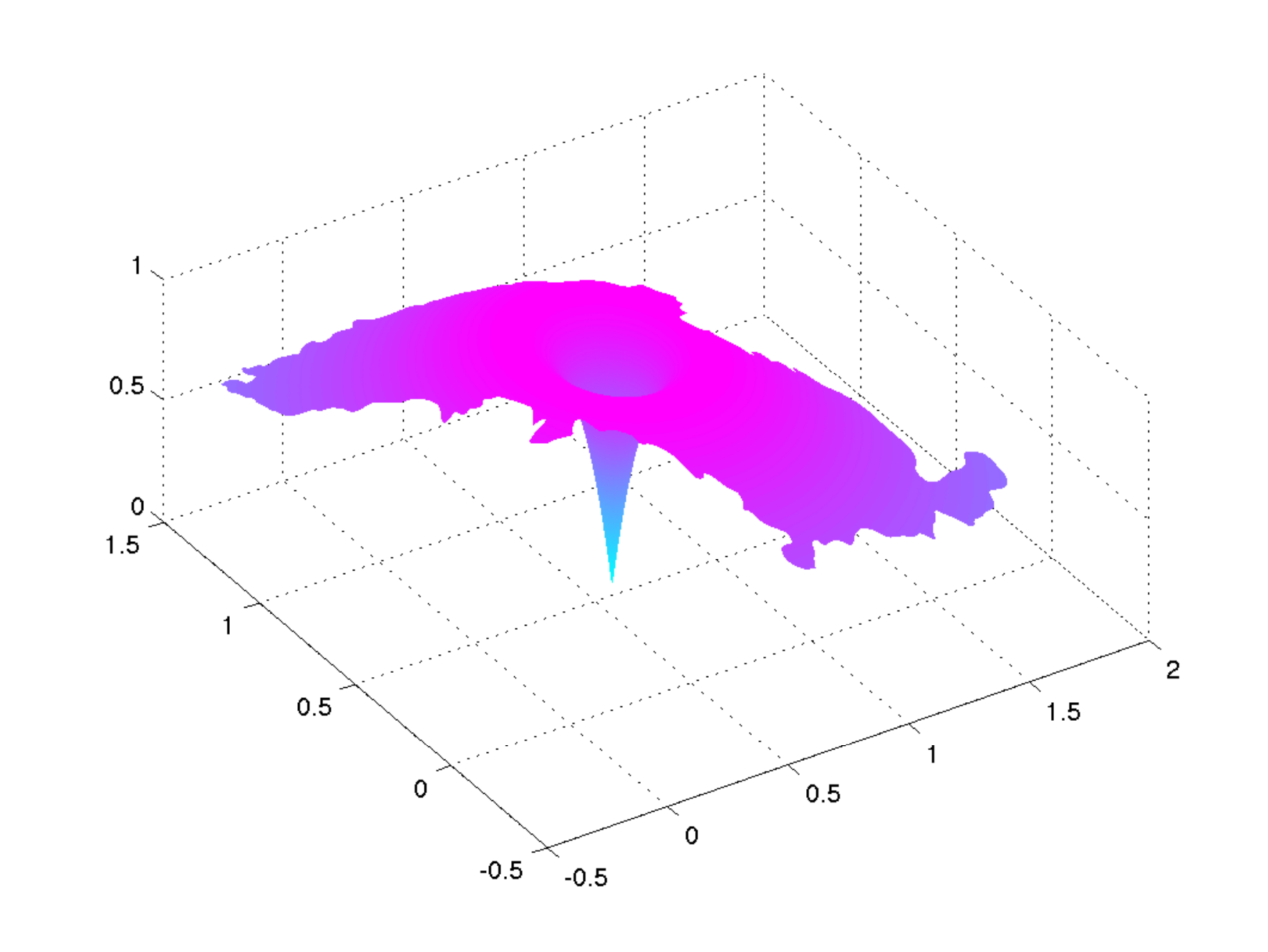}
\includegraphics [scale=0.27]{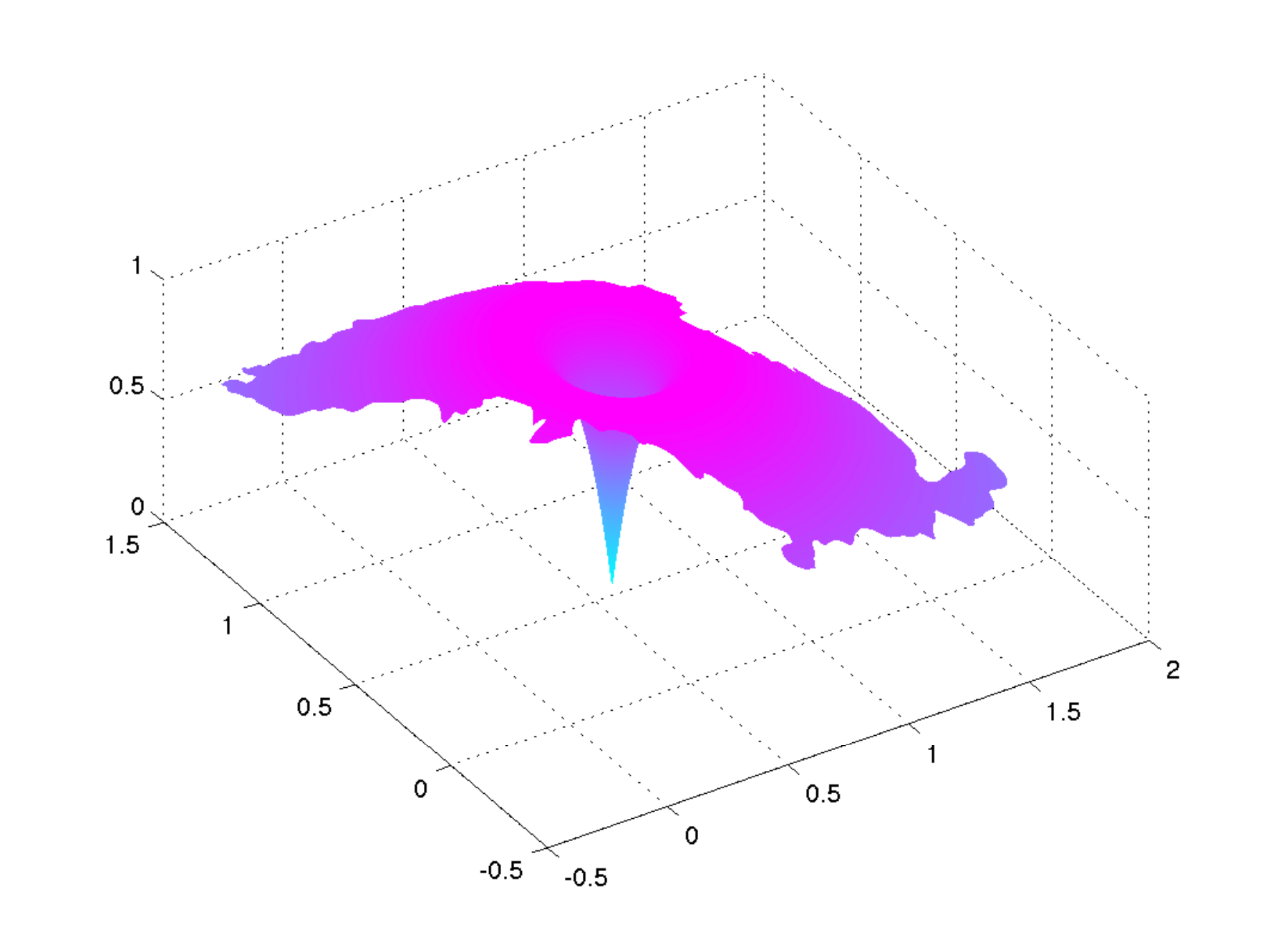}\\
\includegraphics [scale=0.27]{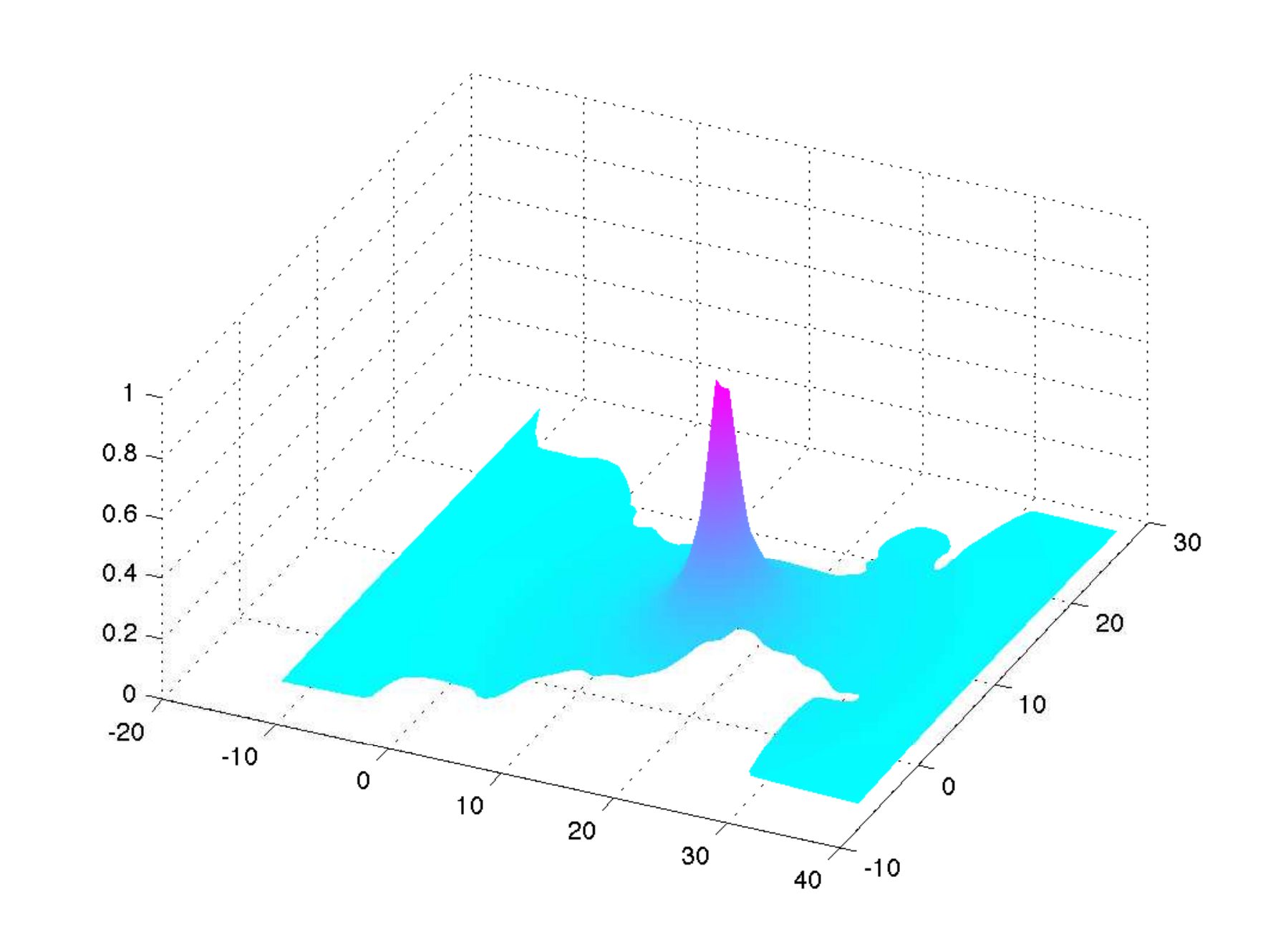}
\includegraphics [scale=0.27]{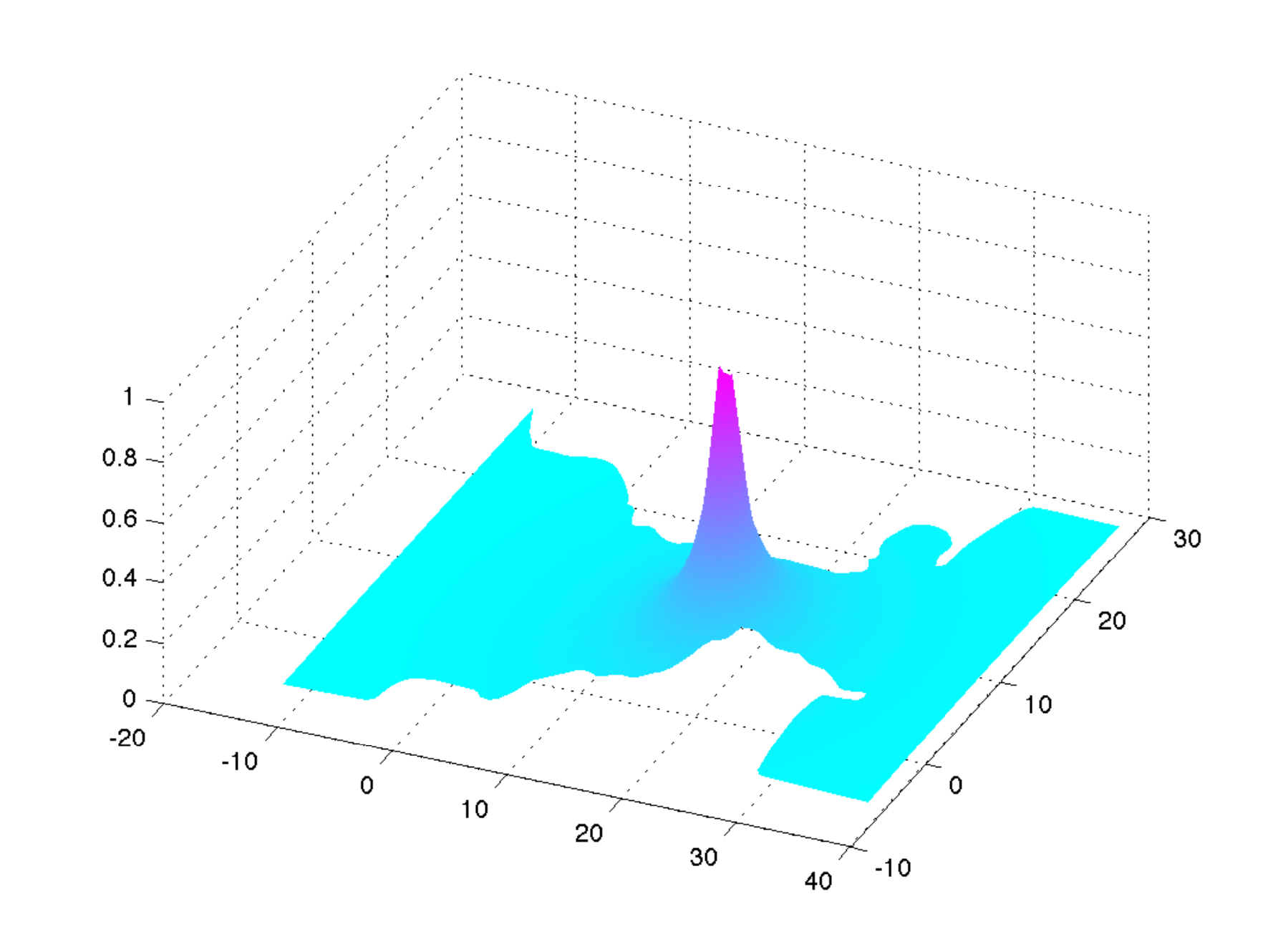}
\end{center}
\caption{ Biquadratic B-spline solution $u^h$ and exact solution (\ref{SolHelm}) of Helmholtz equation with variable frequency. Left column: B-spline approximation $u^h$. Right column: exact solution (\ref{SolHelm}). First row Toba lake, second row Gibraltar channel.}
\label{Fig:SolHemEjBSp-Exact1}
\end{figure}

Initially the sequences of knots $t^{\xi}$ and $t^{\eta}$ are defined by (\ref{tchi_cuad}) and  (\ref{teta_cuad}) respectively. But new knots are inserted depending on the position of the point $(\widetilde{\xi},\widetilde{\eta})$. More precisely, if $ t_i^{\xi} < \widetilde{\xi} <t_{i+1}^{\xi}$ then the knot sequence (\ref{tchi_cuad}) is refined inserting 9 equally spaced knots in the intervals $(t_{i-1}^{\xi},t_i^{\xi}),(t_{i}^{\xi},t_{i+1}^{\xi})$ and $(t_{i+1}^{\xi},t_{i+2}^{\xi})$. Similarly, if $t_j^{\eta} < \widetilde{\eta} <t_{j+1}^{\eta}$ then we insert in (\ref{teta_cuad}) 9 knots equally spaced in the intervals $[t_{j-1}^{\eta},t_j^{\eta}],[t_j^{\eta},t_{j+1}^{\eta}]$ and $[t_{j+1}^{\eta},t_{j+2}^{\eta}]$. If $\widetilde{\xi}$ or $\widetilde{\eta}$ agrees with a knot of the sequences $t^{\xi}$ and $t^{\eta} $ respectively, then we insert 9 equally spaced knots in both intervals of $t^{\xi}$ and $t^{\eta}$ containing the value $\widetilde{\xi}$ and $\widetilde{\eta}$. Finally, since the gradient of the exact solution $u$ is discontinuous in $(x_0,y_0)$ we always insert $\widetilde{\xi}=\widetilde{\eta}=0.5$ as a double knot in $t^{\xi}$ and also as a double knot in $t^{\eta}$.

Figure \ref{Fig:SolHemEjBSp-Exact1} shows the graph of the exact solution $u$ and the approximated biquadratic B-spline solution $u^h$ for two  of the regions reported in table \ref{tablaHelmholtzM1}. For each row, the image in left column is $u^h$ and the image in the right column is $u$. The differences between $u$ and $u^h$ are not appreciable.

In table \ref{tablaHelmholtzM1} we show the errors of the biquadratic B-spline solution $u^h$ for different physical regions. The $L_2$ error oscillates between $10^{-3}$ and $10^{-4}$, but the $H_1$ error is
two orders bigger.

\begin{table}[hbt]
\begin{center}
\begin{tabular}{||c|c|c|c||}
\hline
\hline
Region                   & Degrees of freedom & $L_2\;error$  & $H_1\;error$\\\hline\hline
Havana bay               &  $117 \times 111$         &  0.0038   &   0.2076\\\hline
Toba lake                &  $167 \times 167$         &  0.0014   &   0.0223\\\hline
Gibraltar channel        &  $ 97 \times 157$         &  0.0051   &   0.6613\\\hline
Grijalva channel         &  $119 \times 39 $         &  0.0015   &   0.4784\\\hline
P\'atzcuaro lake         &  $109 \times 109$         & 9.6912e-4 &   0.2712\\\hline
V. de Bravo reservoir    &  $157 \times 157$         & 3.3098e-4 &   0.0152\\
\hline
\hline
\end{tabular}
\caption{\label{tablaHelmholtzM1} Errors of the biquadratic B-spline solution of Helmholtz equation with exact solution (\ref{SolHelm}) on several physical regions.}
\end{center}
\end{table}

\subsubsection{Experiments increasing the number of oscillations of the exact solution. }

As we already mentioned, the parameter $M$ in the expression of $\alpha=\frac{1}{M \pi}$ is the number of oscillations of the exact solution (\ref{SolHelm}). Hence, in order to obtain a good approximation $u^h$ of the exact solution $u$, for values of $M$ greater than 1 we must add more basic functions $B_{i}^3(\xi)$  and $B_{j}^3(\eta)$ different from zero near $\widetilde{\xi}=\widetilde{\eta}=0.5$. In our experiments,
we always insert $\widetilde{\xi}=\widetilde{\eta}=0.5$ as a double knot in $t^{\xi}$ and also as a double knot in $t^{\eta}$. Moreover, a total of $27$ equally spaced knots are inserted in both intervals in $t^{\xi}$ containing $\widetilde{\xi}$. The same procedure is used for inserting knots in $t^{\eta}$.

\begin{figure}[hbt]
\begin{center}
\includegraphics [scale=0.2]{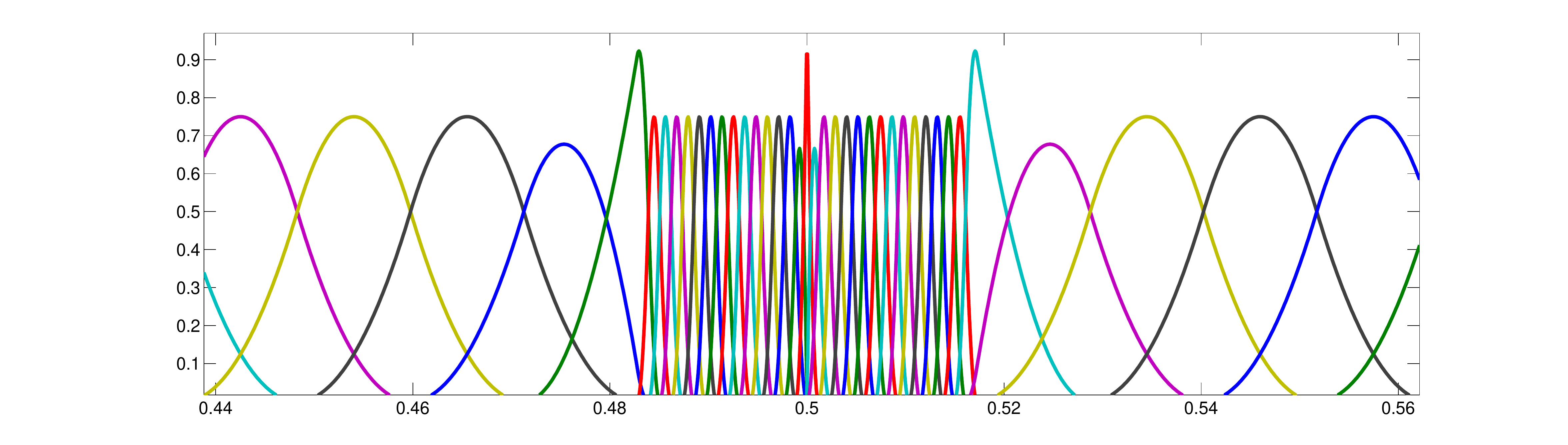}
\label{Fig:baseref}
\end{center}
\caption{Some basis functions $B_{i,t^\xi}^3(\xi)$ after inserted knots near $\xi=0.5$.  }
\label{Fig:basesref}
\end{figure}

In Figure \ref{Fig:basesref} we show some of the basis functions $B_{i}^3(\xi),\;\; 0 \leq \xi \leq 1$. Since the knots are very concentrated in the neighborhood of $\widetilde{\xi}=0.5$, we observe that many basic functions (one for each knot inserted) are different from 0 near this value.

In the next examples we solve the Helmholtz equation with exact solution (\ref{SolHelm}) for $M=2$, $M=3$ and $M=4$. In table \ref{tablaHelmholtzM1234} we show the $L_2$ and the $H_1$ errors of the approximated solution $u^h$, when the physical domain is Havana bay. For comparison, we also include the result of table \ref{tablaHelmholtzM1} for $M=1$. Observe that the number of degrees of freedom, reported in column 2, is bigger for $M>1$ than for $M=1$.

\begin{table}[htb]
\begin{center}
\begin{tabular}{||c|c|c|c||}
\hline
\hline
Number of oscillations ($M$) & Degrees of freedom &  $L_2\;error$ &  $H_1\;error$ \\\hline\hline
1 &  $117 \times 111$  &  0.0038  &   0.2076\\\hline
2 &  $171 \times 165$  &  0.0018  &   0.3394\\\hline
3 &  $171 \times 165$  &  0.0038  &   1.1035\\\hline
4 &  $171 \times 165$  &  0.0069  &   1.987\\
\hline
\hline
\end{tabular}
\caption{\label{tablaHelmholtzM1234} Errors on Havana bay of the biquadratic B-spline solution of Helmholtz equation, with exact solution (\ref{SolHelm}) for increasing number $M$ of oscillations.}
\end{center}
\end{table}

\begin{figure}[hbt]
\begin{center}
\includegraphics [scale=0.23]{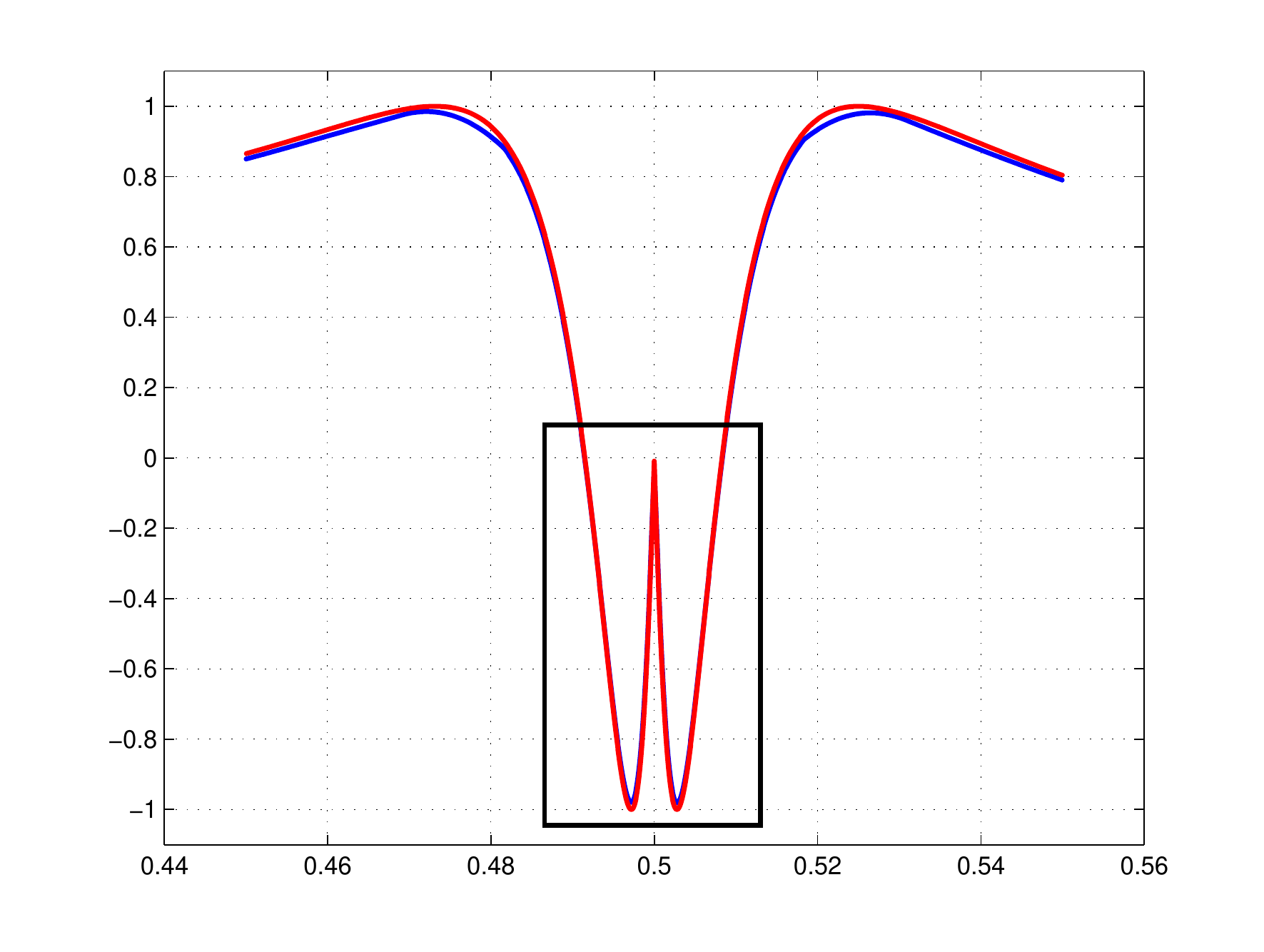}
\includegraphics [scale=0.23]{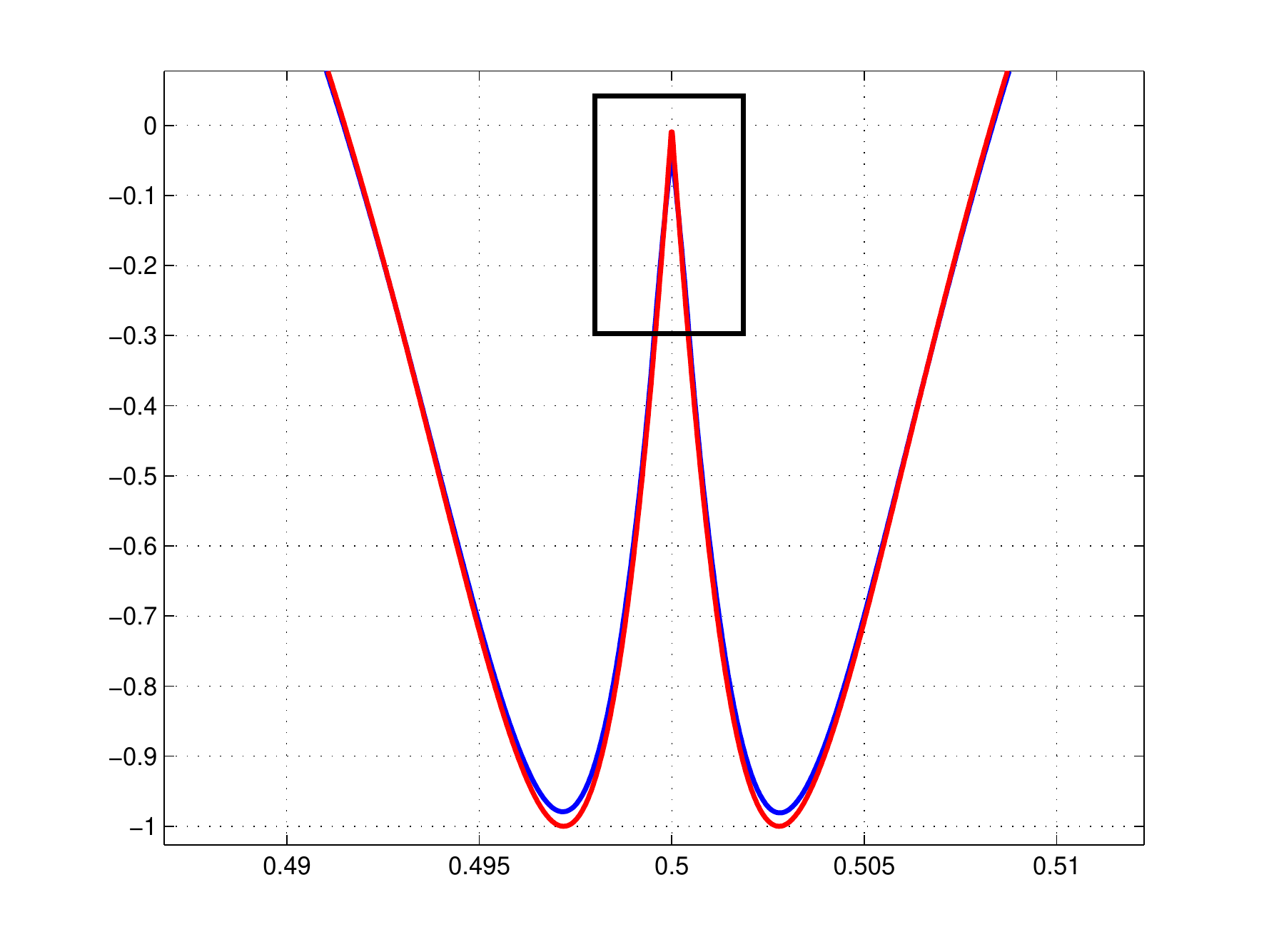}
\includegraphics [scale=0.23]{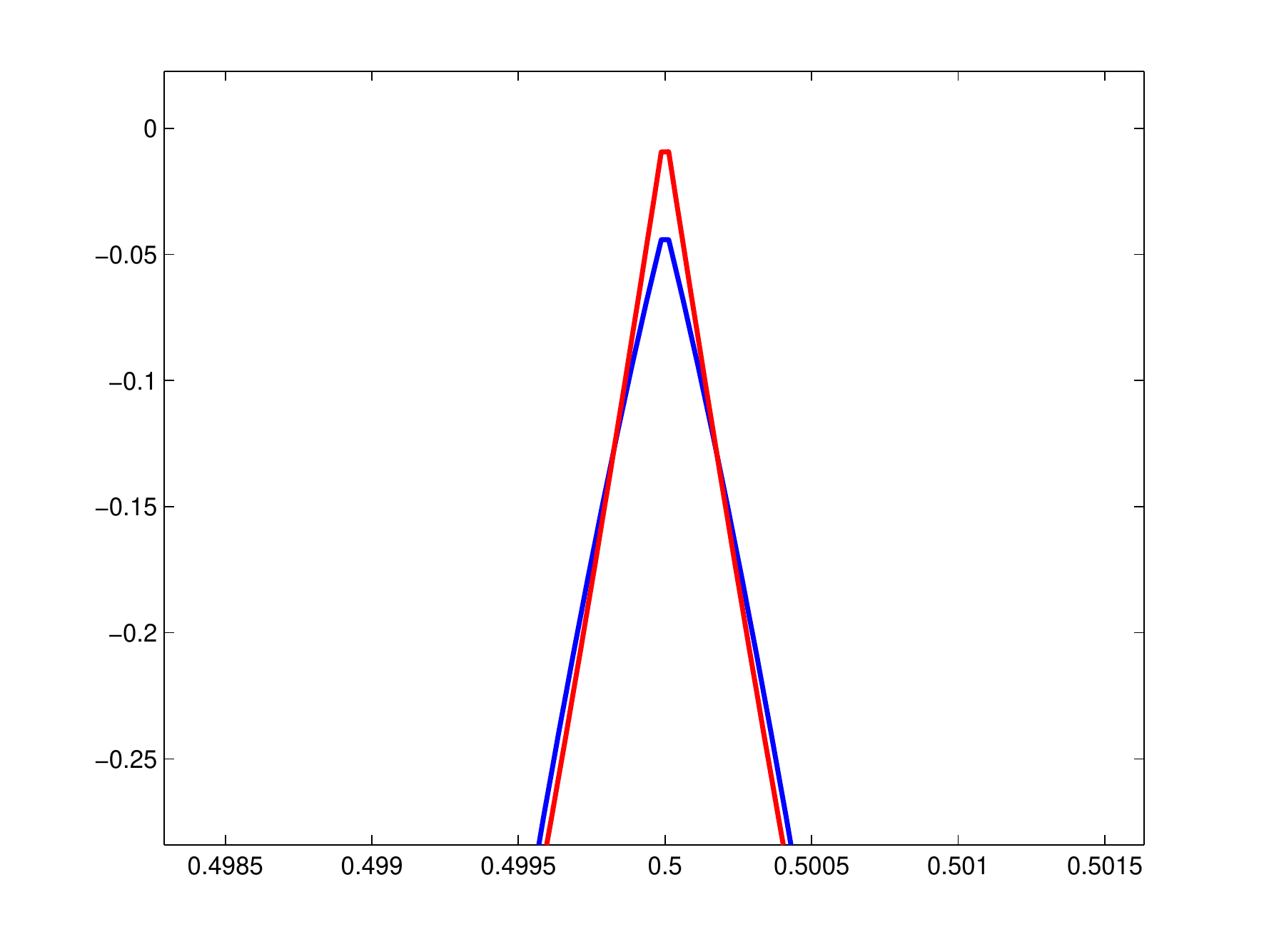}
\includegraphics [scale=0.23]{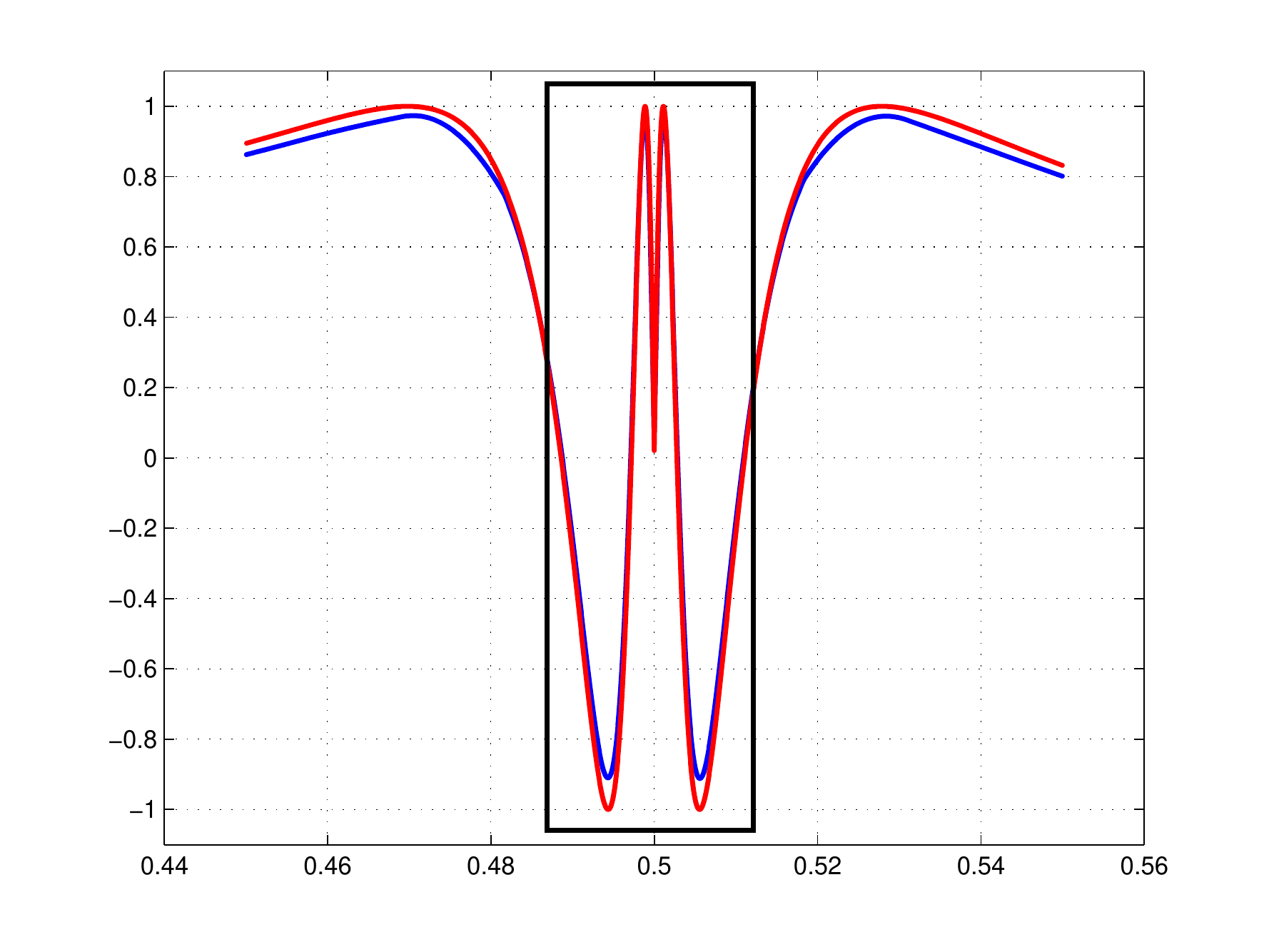}
\includegraphics [scale=0.23]{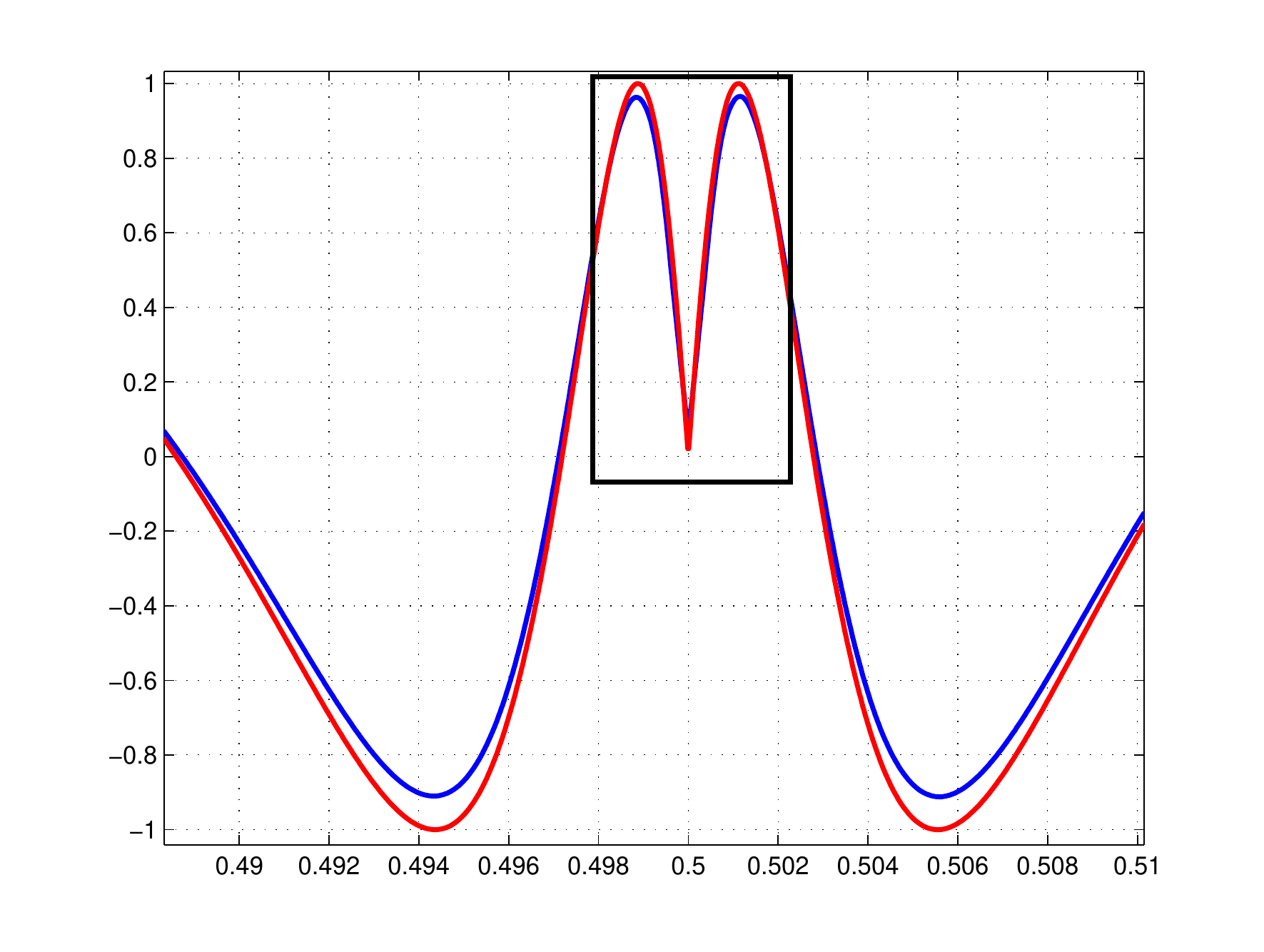}
\includegraphics [scale=0.23]{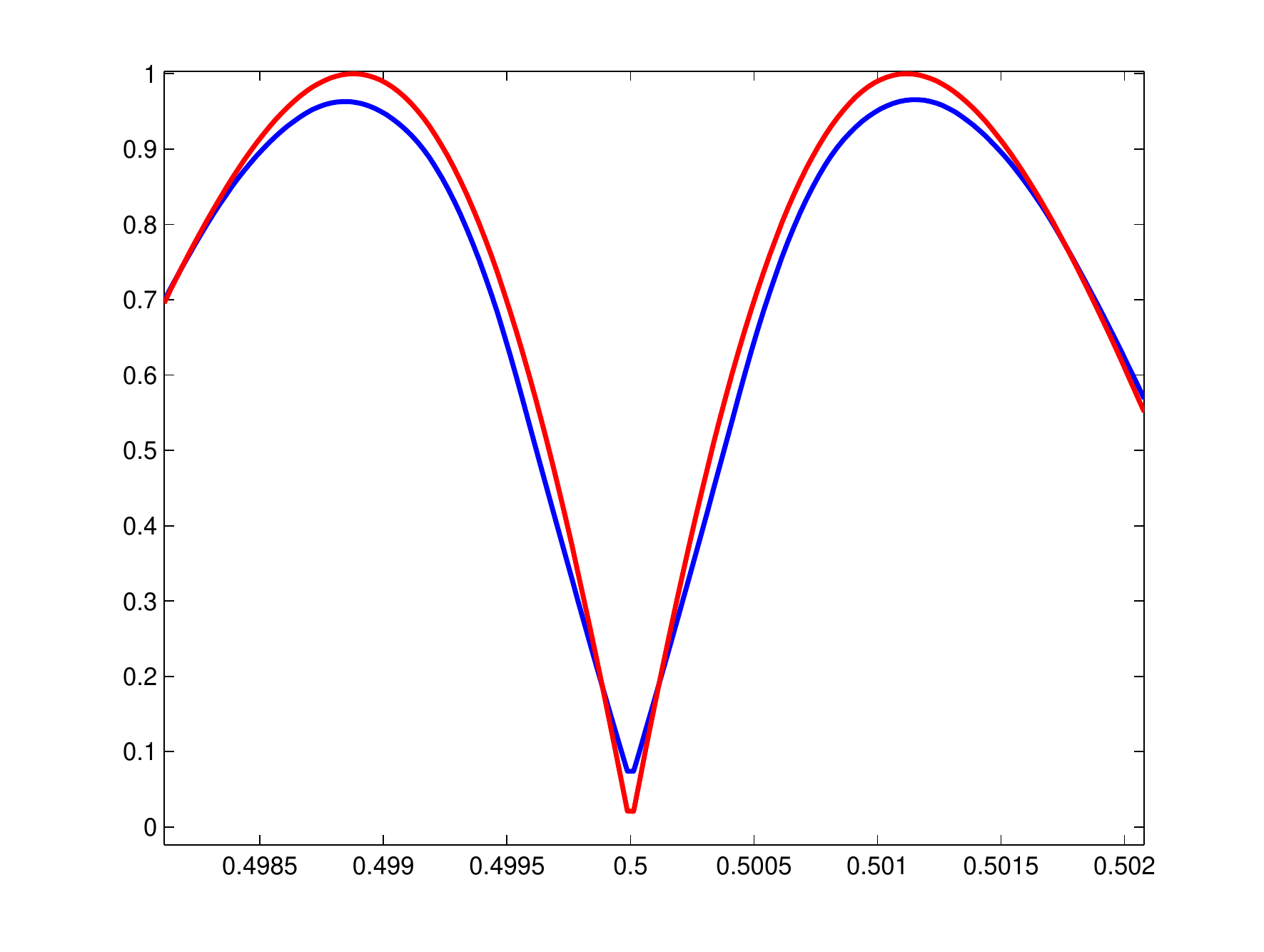}
\includegraphics [scale=0.23]{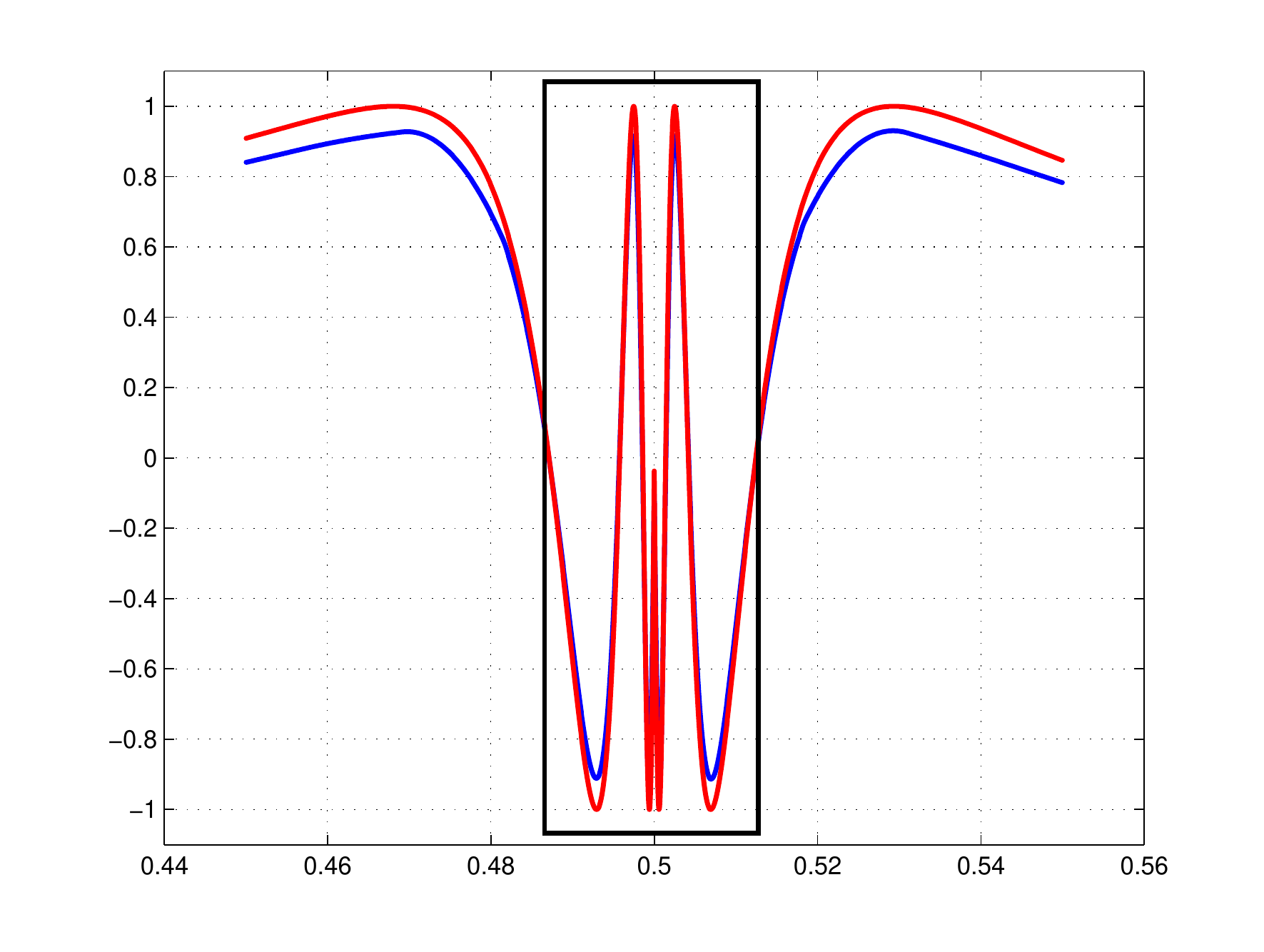}
\includegraphics [scale=0.23]{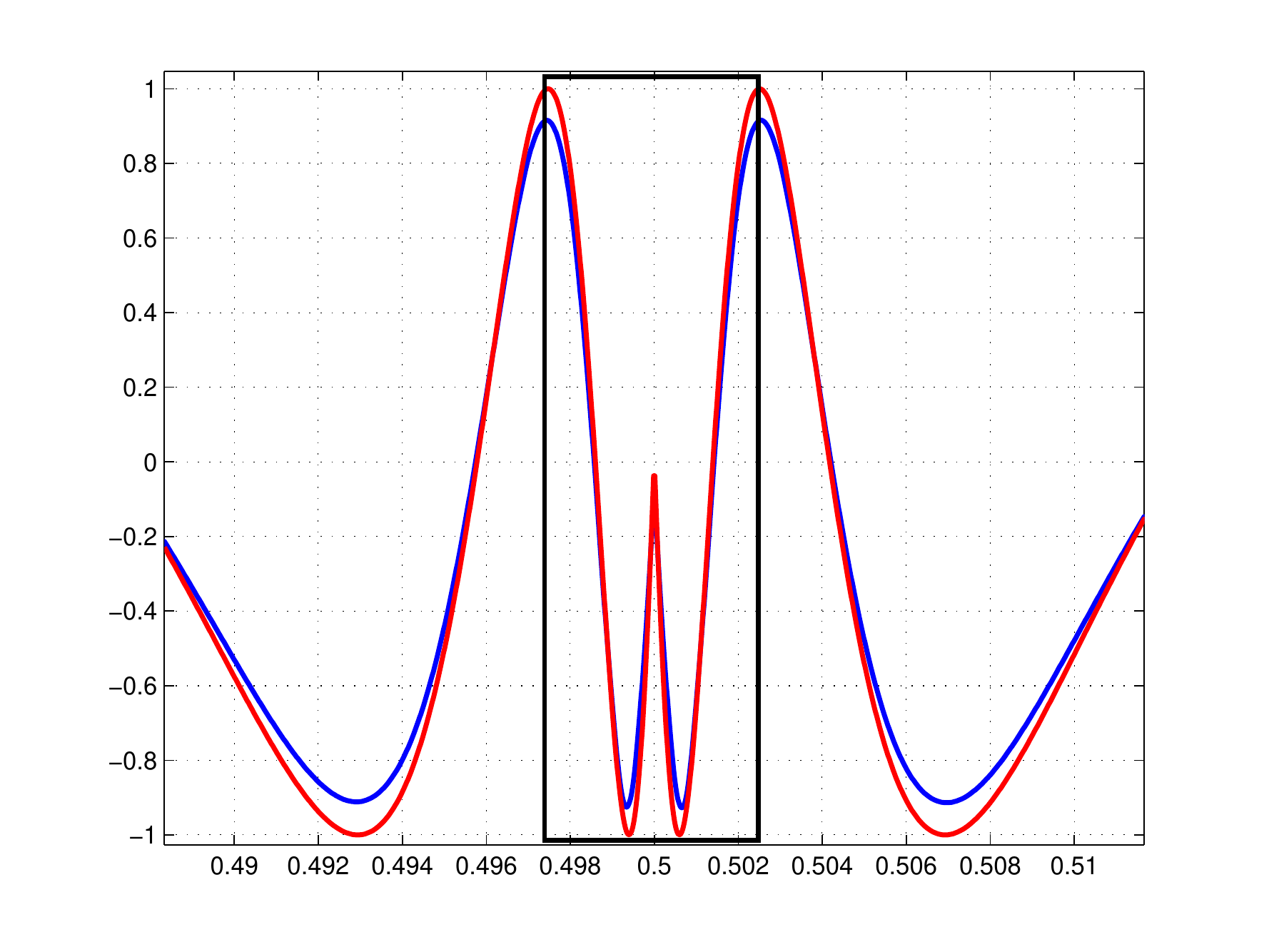}
\includegraphics [scale=0.23]{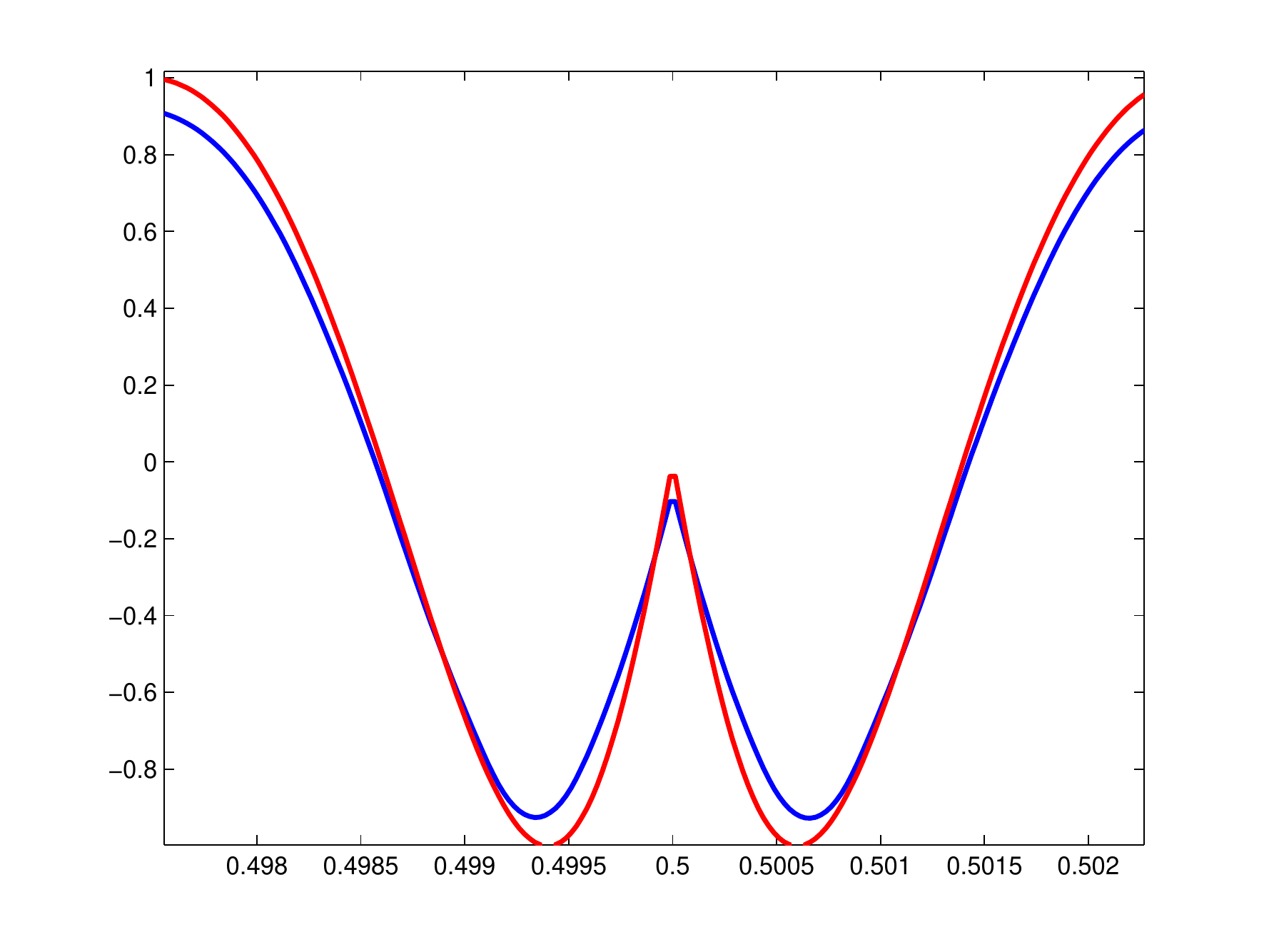}
\end{center}
\caption{Restriction to a parametric line of the exact (red) and approximated (blue) solutions of the Helmholtz equation with variable frequency on Havana bay. Each rectangle area is amplified on right image in the same row. First row $M=2$, second row $M=3$ and third row $M=4$.}
\label{Fig:BahHabSolHelmOscM234}
\end{figure}
Figure \ref{Fig:BahHabSolHelmOscM234} shows the graph of the functions $u(\mathbf{F}(\xi,\eta))$ and $u^h(\mathbf{F}(\xi,\eta))$ for Havana bay, both restricted to the parametric line $\xi=\eta$, which contains the pre-image of the singular point $(x_0,y_0)$. More precisely, the red graph shows the curve $u(\mathbf{F}(\xi,\xi))$, while the blue graph represents the function $u^h(\mathbf{F}(\xi,\xi))$. The first row corresponds to the solution with $M=2$ oscillations, the second and third rows correspond to $M=3$ and $M=4$ respectively. In each row, the black rectangle area in the graph is zoomed on right image. Observe that the oscillations are located is a very narrow segment. Moreover, the approximated solution $u^h$ reproduces the behavior of $u$ including the non differentiability in the point $\widetilde{\xi}=0.5$.

\subsubsection{Convergence study}

To study the convergence of the isogeometric approach we solve in this section the Helmholtz equation, where $k(x,y)$ and $f(x,y)$ are given by (\ref{kSho}) and (\ref{fxy}) respectively, $\Omega$ is a Jigsaw puzzle region given in \cite{Grav14} and $\alpha=\frac{1}{\pi}$. In table \ref{tablaconv} we report the $L_2$ and the $H_1$ errors for increasing values of the number $N = n \times m$ of degrees of freedom. The $i$-th row of table \ref{tablaconv} corresponds to a sequence of knots which is obtained refining uniformly $i-1$ times the initial uniform sequences (\ref{tchi_cuad}) and (\ref{teta_cuad}) and introducing later equally spaced knots in the intervals containing the parametric value $0.5$ in each direction (see the previous section). As we observe the $L_2$ and the $H_1$ errors decrease as the number of degrees of freedom increases.

\begin{table}[htb]
\begin{center}
\begin{tabular}{||c|c|c||}
\hline
\hline
Degrees of freedom $n \times m$           &  $L_2\;error$ &  $H_1\;error$ \\\hline\hline
        $ 13  \times 13$    & 0.2024     & 3.4531  \\\hline
        $ 25  \times 25$    & 0.0558     & 1.6077  \\\hline
        $ 45  \times 45$    & 0.0168     & 0.7877  \\\hline
        $ 81  \times 81$    & 0.0061     & 0.1234  \\\hline
        $ 149 \times 149$   & 0.0024     & 0.1213  \\\hline
        $ 281 \times 281$   & 0.0010     & 0.0404  \\\hline
        $ 537 \times 537$   & 0.0005     & 0.0301  \\
\hline
\hline
\end{tabular}
\caption{\label{tablaconv} Errors on Jigsaw puzzle region of the biquadratic B-spline solution of Helmholtz equation, with exact solution (\ref{SolHelm}) for increasing number of degrees of freedom.}
\end{center}
\end{table}

%

\begin{figure}[htb]
\center
\includegraphics[scale=0.22]{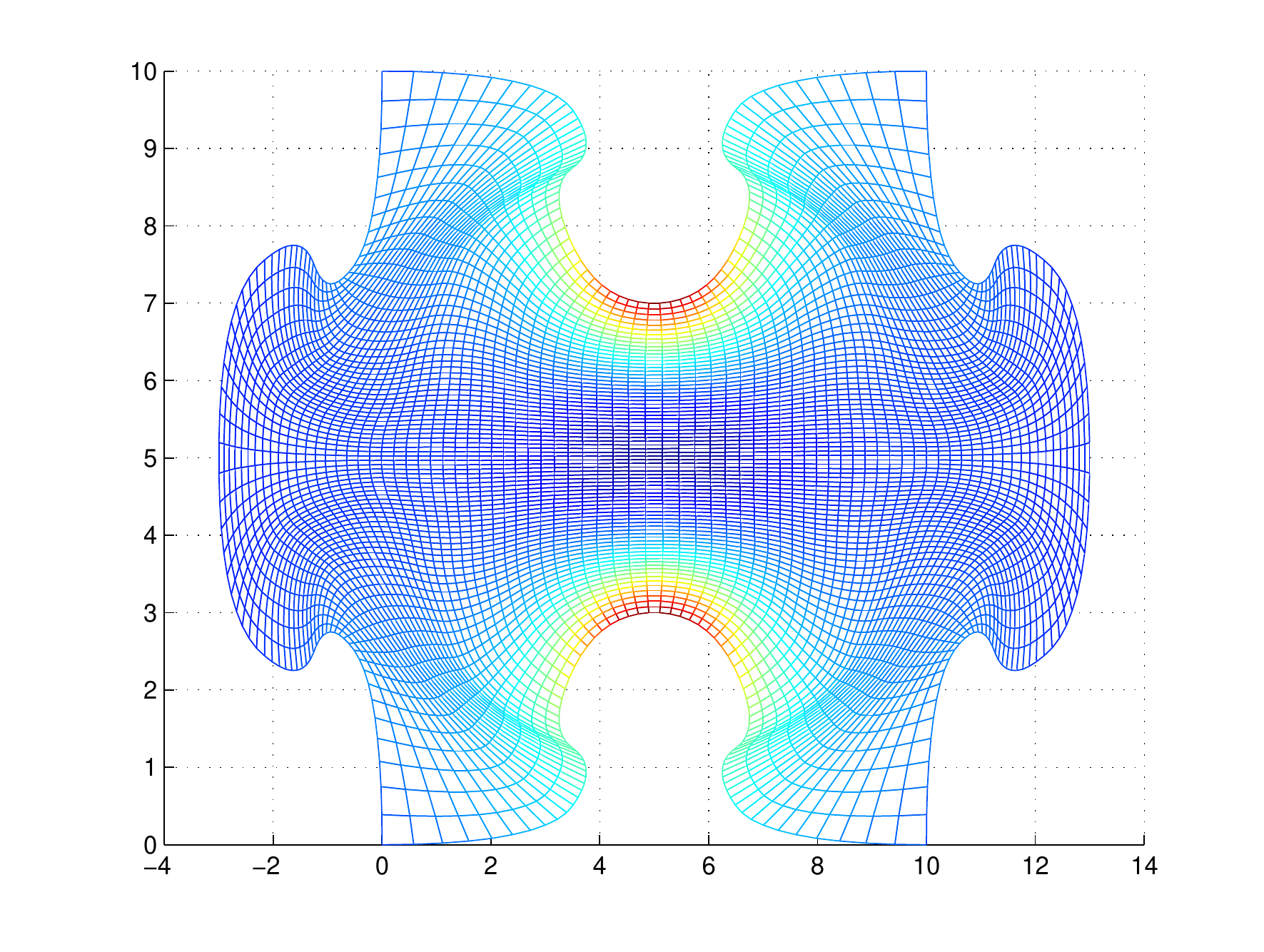}
\includegraphics[scale=0.22]{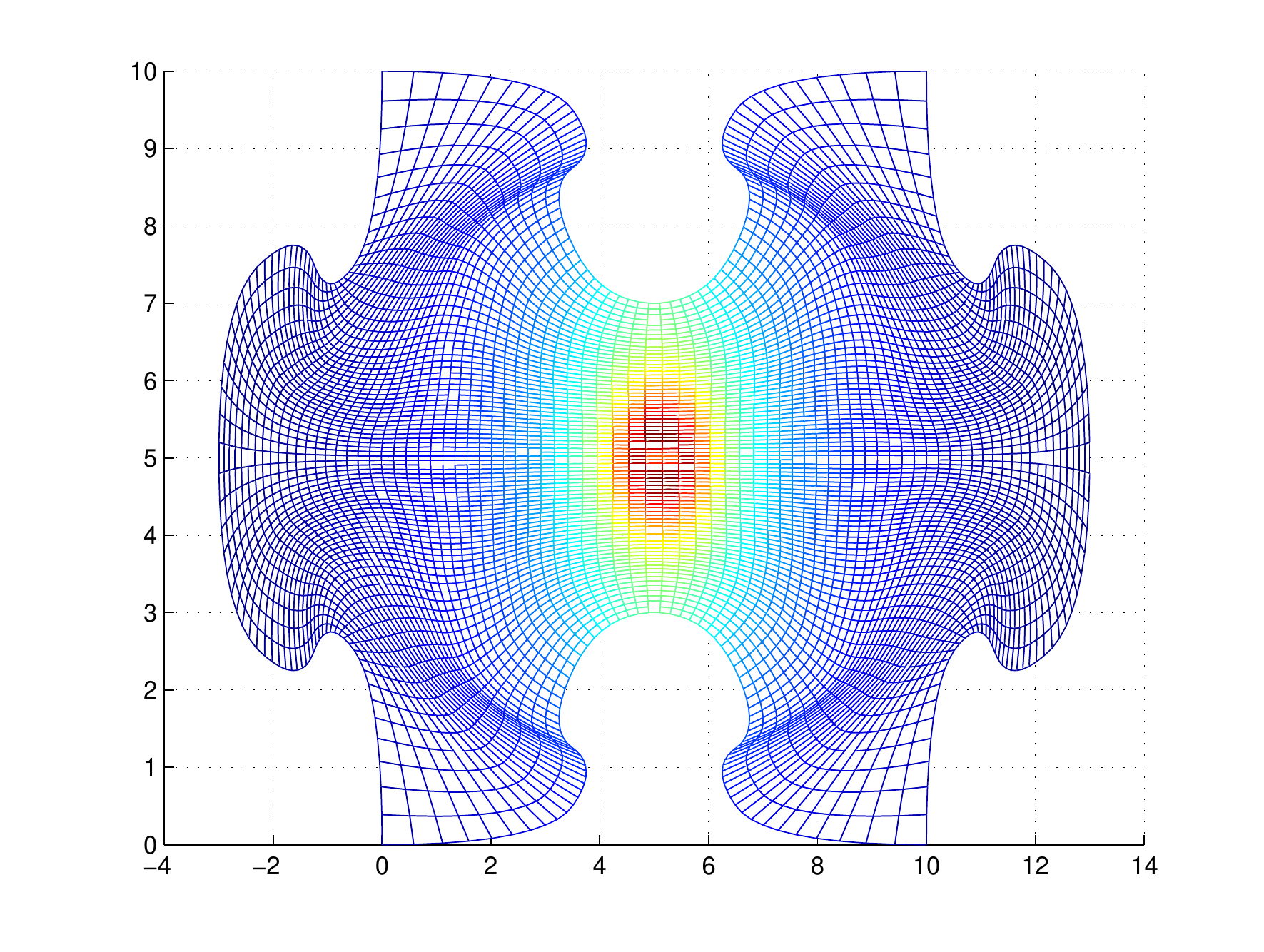}
\includegraphics[scale=0.22]{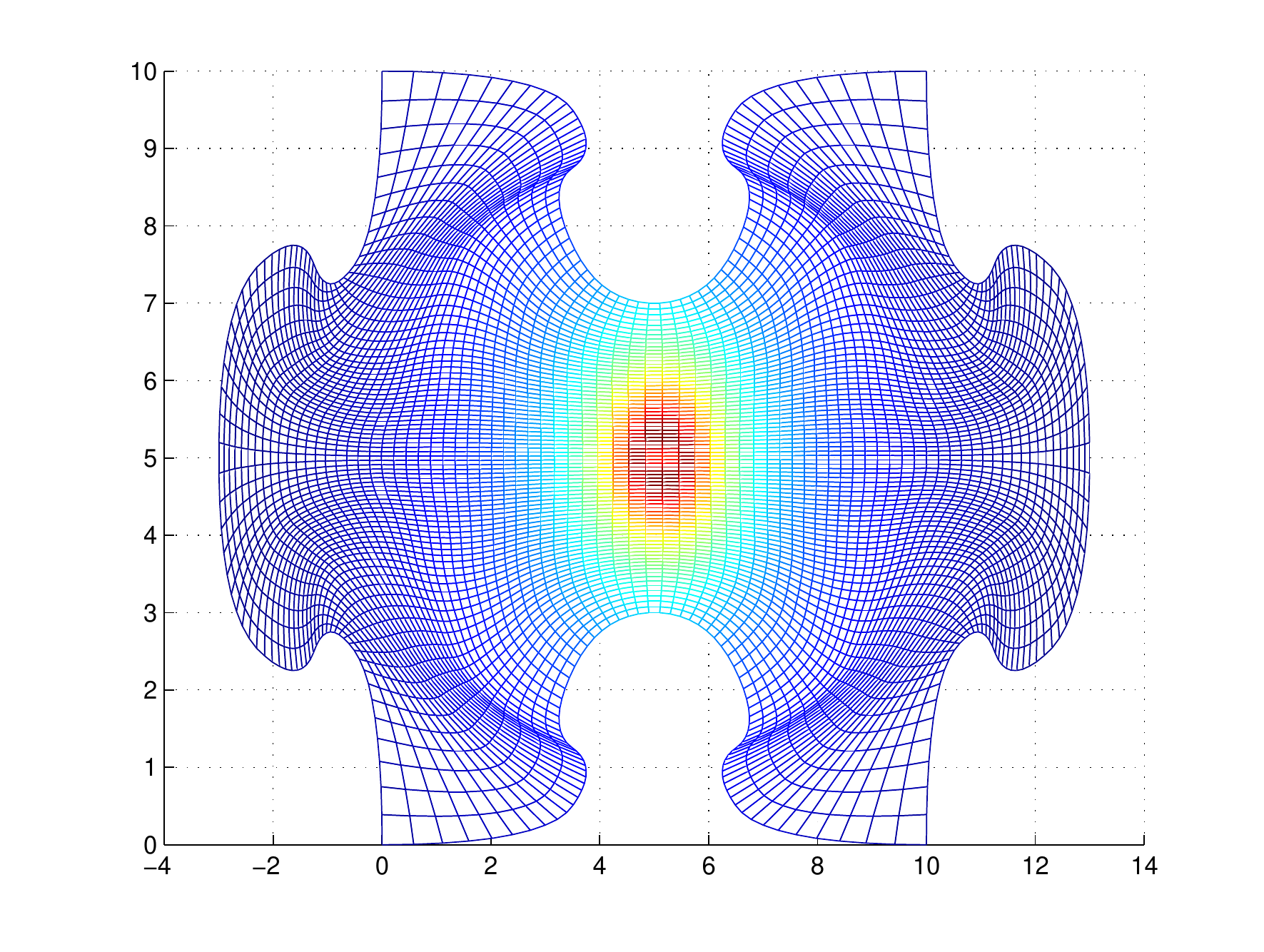}
\caption{2D views of the approximated solution of the Helmholtz equation when we increase the number of degrees of freedom. From left to right the total number $N$ of degrees of freedom is: $169,\; 625$ and $288369$.}
\label{Fig:Puzzle}
\end{figure}

In Figure \ref{Fig:Puzzle} we show a 2D view of the approximated solution for three of the cases reported in table \ref{tablaconv}. There are almost no differences between the approximated solution with $537 \times 537=288369$ degrees of freedom and the exact solution.

\section{Conclusions}

The solution of partial differential equations using the IgA approach has several advantages in comparison with the classical finite element method. One of them is that the boundary of the physical domain is represented exactly. This is specially important when the domain is a region with irregular boundary and high error would be introduced if the boundary has to be approximated. On the other hand, IgA approach is able to produce smoother solutions having at the same time some singular points and high oscillations.

In this paper we have focussed our attention on the solution of Helmholtz equation with Dirichlet boundary condition. To approximate the solution we have used biquadratic B-spline functions selecting carefully the number and position of the knots, in such away that the approximated solution reproduces the behavior of the exact solution, even when the last one has singular points or zones of high oscillations. In this sense, our numerical experiences show that IgA approach can be successfully used to solve difficult cases of the Helmholtz equation. The success of the method also depends on the quality of the parametrization of the physical domain, specially when it has very irregular boundary.

As a future work we plan to solve the Helmholtz equation using IgA approach, when the frequency is a big positive constant representing the wave number. This problem, very important in acoustics and other applications, is difficult to solve with classical FEM. Our intention is to show that IgA approach is a better option to obtain good approximated solutions.


\end{document}